\documentclass[a4paper,reqno,11pt]{amsart}
\usepackage{graphicx,color,marginnote}
\usepackage{footmisc}
\usepackage{amsmath}

\usepackage{amsfonts}
\usepackage{amssymb}
\usepackage{bbm}
\usepackage{epstopdf}
\usepackage{color}
\numberwithin{equation}{section}
\usepackage[marginparwidth=3cm, marginparsep=.5cm]{geometry}

\newcommand{\bbR}{\mathbb{R}}

\usepackage{amsthm}
\newtheorem{Theorem}{Theorem}

\newtheorem{Remark}[Theorem]{Remark}

\newtheorem{Assumption}[Theorem]{Assumption}

\newtheorem{Definition}[Theorem]{Definition}
\newtheorem{Proposition}[Theorem]{Proposition}

\title[Hydrodynamic limit for dimer dynamics]{Hydrodynamic limit
  equation for a lozenge tiling Glauber dynamics}
\author{Beno\^it  Laslier}
\address{Statslab, 
Centre for Mathematical Sciences, 
Wilberforce Road, 
Cambridge CB3 0WA - UK\\
and\\
LPMA - Univ. Paris Diderot,
Bâtiment Sophie Germain,
avenue de France,
75013 Paris, France
\\
E-mail: laslier@math.univ-paris-diderot.fr}
\author{Fabio Lucio Toninelli}
\address{Univ Lyon 1, Universit\'e Claude Bernard Lyon 1, CNRS UMR 5208, Institut Camille Jordan, 43 blvd. du 11 novembre 1918, F-69622 Villeurbanne cedex, France
 \\
E-mail: toninelli@math.univ-lyon1.fr}

\begin{document}

\maketitle

\begin{abstract}
  We study a reversible continuous-time Markov dynamics on lozenge
  tilings of the plane, introduced by Luby et al. \cite{LRS}.  Single
  updates consist in concatenations of $n$ elementary lozenge
  rotations at adjacent vertices. The dynamics can also be seen as a
  reversible stochastic interface evolution. When the update rate is
  chosen proportional to $1/n$, the dynamics is known to enjoy
  especially nice features: a certain Hamming distance between
  configurations contracts with time on average \cite{LRS} and the
  relaxation time of the Markov chain is  diffusive
  \cite{Wilson}, growing like the square of the diameter of the
  system. Here, we present another remarkable feature of this
  dynamics, namely we derive, in the diffusive time scale, a fully
  explicit hydrodynamic limit equation for the height function (in the
  form of a non-linear parabolic PDE). While this equation
  \emph{cannot} be written as a gradient flow w.r.t. a surface energy
  functional, it has nice analytic properties, for instance it
  contracts the $\mathbb L^2$ distance between solutions. The mobility
  coefficient $\mu$ in the equation has non-trivial but explicit
  dependence on the interface slope and, interestingly, is directly
  related to the system's surface free energy.  The derivation of the
  hydrodynamic limit is not fully rigorous, in that it relies on an
  unproven assumption of local equilibrium.

\end{abstract}

\section{Introduction}

The large-scale  time evolution of interfaces separating different thermodynamic phases
is a classical subject in statistical mechanics. A first natural goal
is that of obtaining a hydrodynamic limit \cite{KL,SpohnLibro,Fu}: take an initial interface
configuration that approximates a macroscopic smooth profile, let it
evolve via a microscopic Markovian Glauber-type dynamics  that, at the lattice level,
follows simple local rules and, rescaling time and space properly,
prove that the interface converges to the solution of a deterministic
PDE. If the two thermodynamic phases separated by the interface are at
coexistence, i.e. if they have the same bulk free energy, we expect
the correct time rescaling to be diffusive and the limit equation to be a
parabolic PDE, in general a non-linear one, of the form
\begin{eqnarray}
  \label{eq:eqgenerale}
  \partial_t \phi(x,t)=-\mu(\nabla \phi(x,t))\frac{\delta F[\phi]}{\delta \phi(x,t)}.
\end{eqnarray}
Here, $F[\phi]$ is the surface tension functional, that is a purely
equilibrium quantity, while $\mu(\nabla \phi)>0$ is the (in general slope-dependent) interface
mobility coefficient, that depends on the generator of the Markov chain.

To obtain a mathematically simpler model, the interface is often
described at the microscopic level by a $d+1$-dimensional height
function (``effective interface''), i.e. the graph of a function from
$\mathbb Z^d$ to $\mathbb R$ (or to $\mathbb Z$ in the case of
discrete interface models). Here, $d+1$ is the dimension of space
where the thermodynamic system of interest lives and of course the
physically most relevant case is $d=2$.  In the
``effective interface'' approximation, the internal structure of the
two bulk phases is forgotten and the occurrence of 
 interface overhangs is
entirely neglected. Despite this somewhat drastic simplification, and
despite the fact that the phenomenological picture behind the expected
hydrodynamic limit is rather clear \cite{Spohn}, most effective
interface dynamics remain mathematically intractable and rigorous
progress is very limited, especially for $ d>1$. One notable exception
is that of the Langevin dynamics for the Ginzburg-Landau model with
symmetric and strictly convex potential, where a rigorous derivation
of the hydrodynamic limit was obtained for any $d\ge1$ by Funaki and
Spohn \cite{FuSpo} (see also \cite{Nishi} who extended \cite{FuSpo}
beyond the case of periodic boundary conditions).

Leaving aside the problem of rigorously proving the hydrodynamic
limit, even the more modest goal of guessing the exact form of the
limit PDE is in general out of reach, except for lucky exceptions (the
Ginzburg-Landau model being one of them) where the dynamics satisfies
some form of ``gradient condition'' \cite{SpohnLibro,KL}
which allows to obtain a simple formula for the interface mobility
$\mu(\cdot)$, involving only equal-time equilibrium averages.

\medskip

The goal of the present work is to present a Markov chain for a  discrete
interface model in dimension $(2+1)$ and to show that it should admit a hydrodynamic limit
that is fully explicit and non-trivial ($\mu(\nabla\phi)$ is a non-linear function of the interface slope). 
We comment below on what is missing in order to turn our arguments into a rigorous proof.

\begin{figure}[h]
  \includegraphics[width=7cm]{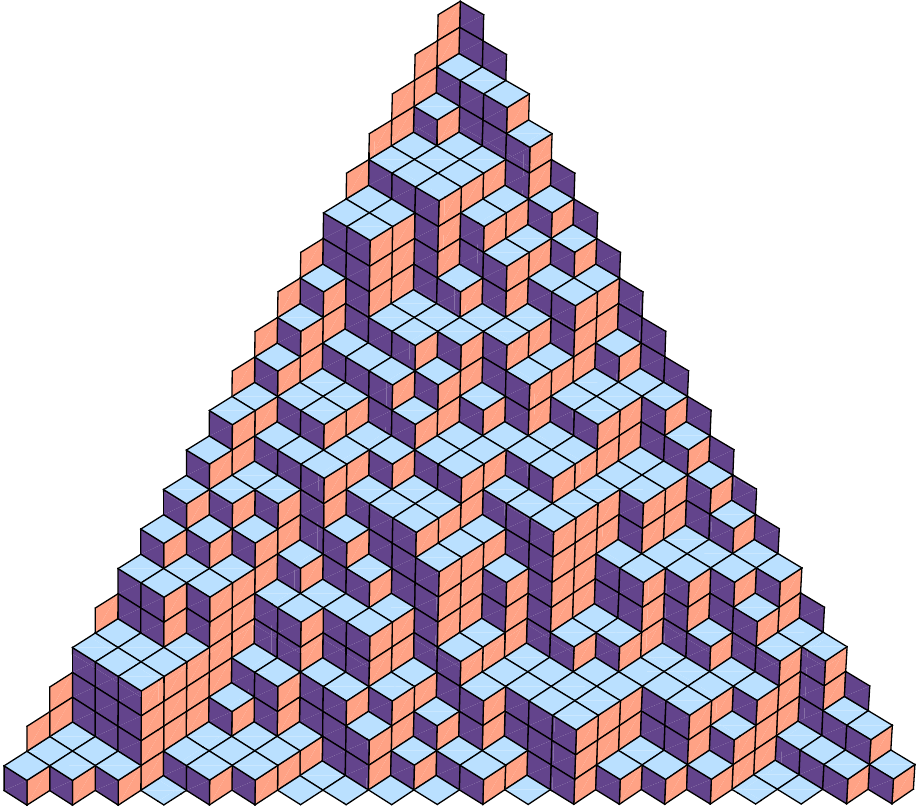}
\caption{A (portion of) lozenge tiling. Viewing the tiling as a stacking of
  cubes standing on the horizontal $(x,y)$ plane,
  note that the height of columns (each with a blue lozenge on top) is
weakly decreasing in both $x$ and $y$ directions.}
\label{fig:cubi}
\end{figure}

Before introducing the interface dynamics we are interested in, we
 make a brief detour to motivate the reader.  A class of
discrete interface dynamics that attracted much attention lately are Glauber
dynamics of dimer models, in particular lozenge tilings of the plane
\cite{LRS,Wilson,CMST,CMT,LTCMP}. Such tilings are in
bijection with $(2+1)$-dimensional discrete surfaces obtained as
a monotone stacking of elementary cubes in $\mathbb R^3$, see Figure
\ref{fig:cubi}.  
Here, ``monotone'' means that the heights of columns of cubes, indexed
by the coordinates $(x,y)$ of their orthogonal projection on the horizontal
plane, 
are weakly decreasing both w.r.t. $x$ and $y$. 

\begin{figure}[h]
  \includegraphics[width=5cm]{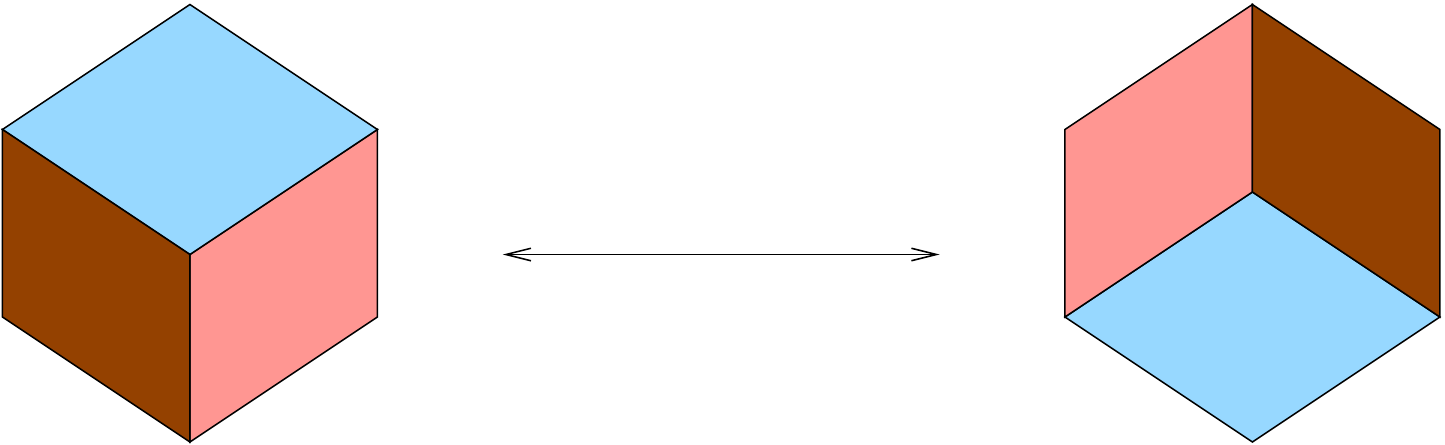}
\caption{The updates of the single-flip dynamics. Both have transition
rate $1$.}
\label{fig:flip}
\end{figure}

The most natural
reversible Markov dynamics on such tilings is the one whose elementary
moves are rotations by an angle $\pm 60^\circ$ (with transition rate $1$) of
three lozenges sharing a common vertex, see Figure \ref{fig:flip}.  
This will
be called the ``single-flip dynamics'' in the following.
As discussed for instance in \cite{CMST,CMT}, the single-flip dynamics
coincides with the zero-temperature Glauber dynamics of $+/-$ spin
interfaces of the three-dimensional Ising model with zero magnetic
field, where spins flip one by one. In terms of monotone stacking of
cubes, the dynamics corresponds to adding/removing a cube to/from 
a column, with transition rate $1$,
provided the cube stacking remains monotone after the update.
Recently it was proven that, if we restrict the single-flip dynamics to domains of
diameter $L$, under certain restrictions on the domain shape the
mixing time is of order $O(L^{2+o(1)})$ as
$L\to\infty$ \cite{CMT,LTCMP}.  These results support the idea that the
correct time-scale to observe a hydrodynamic limit should be diffusive (i.e. that we should
rescale time by $L^2$ to see a macroscopic evolution) but they are far
from being sufficient for proving convergence to a limit PDE.

In the present work, we study two modifications of the single-flip
dynamics, where one allows a number $n\ge1$ of cubes to be added/removed
from a column in each move, again subject to the constraint that the
update is legal (i.e. that the resulting configuration is still a monotone
stacking of cubes). If the rates are carefully chosen as functions of
$n$, the dynamics enjoys much nicer properties than the single-flip
one. The first dynamics we will consider is the one where the transition
rate of a legal update involving the addition/removal of $n$ cubes is
proportional to $1/n$; in the second dynamics, instead, with rate $1$
the height of each column of cubes is resampled from the uniform
distribution on all the allowed values it can take given the height of
neighboring columns.  See Definitions \ref{def1} and \ref{def2} below
for more details.  The former dynamics was originally introduced  in
\cite{LRS}, and the latter in \cite{CMST}. Both are known to satisfy
the special property that the volume difference between two configuration is (on
average) decreasing with time, which allows to deduce that the mixing
time is at most polynomial in $L$ \cite{LRS}. Moreover, it was proven
in \cite{Wilson} that the inverse spectral gap of the dynamics is
$O(L^2)$ and that, in special domains, a certain one-dimensional
projection of the height function satisfies on average the
one-dimensional discrete heat equation.

Here we show that, under a reasonable but unproven assumption of local
equilibrium, one can obtain the explicit form of the hydrodynamic
limit equation for the height function, see Eqs. \eqref{eq:PDE},
\eqref{eq:4} and \eqref{eq:PDEter} below. Actually, the hydrodynamic
equation turns out to be the same for both variants of the
dynamics. Obtaining such an explicit expression for the hydrodynamic
equation is a somewhat surprising fact; indeed, let us stress that in
general (for instance, for the single-flip dynamics) the assumption of
local equilibrium is \emph{not} sufficient to guess the limit
equation: knowledge of \emph{corrections} to local equilibrium is also
necessary. There is a general heuristic formula \cite{Spohn} for the
mobility coefficient $\mu(\nabla \phi)$ which is a variant of
Green-Kubo formula. It is given as the sum of two terms, one involving
only local averages in the stationary state of slope $\nabla \phi$ and
the second involving a time-integral of time-space correlations in the
stationary state. The latter term cannot in general be computed as it
would require a closed form for  space-time correlations. However in
lucky cases (like ours, see
Section \ref{sec:GK}) this term happens to be zero due to a summation by
parts at the discrete level.

As we already mentioned, our derivation of the hydrodynamic limit
relies on an unproven assumption of local equilibrium. There are
various difficulties in proving such assumption, and the direct
application of standard entropy techniques (see e.g. \cite{KL}) seems
out of question, in particular because the stationary measures of the
model exhibit long-range correlations. The adaptation of the so-called
$H^{-1}$ method employed in \cite{FuSpo,Nishi} looks also challenging:
technically a non-trivial difficulty is to get some a-priori 
control of interface gradients during the evolution (see Remark \ref{rem:tec} below for more details).  In
\cite{FuSpo,Nishi} an important role in this respect was played by
strict convexity of the potential, that fails in our case.  
However, in the case where the system has periodic boundary
conditions, in a forthcoming work \cite{cf:hydro} we manage to
overcome these difficulties and to prove rigorously the validity of
the hydrodynamic limit.

The hydrodynamic equation has nice analytic features. While it is
\emph{not} in the form of the gradient flow with respect of a surface
free energy functional, it can be written in a divergence form
(cf. \eqref{eq:PDEdiv}) that allows to show (see Section \ref{sec:l2})
that the $\mathbb L^2$ distance between solutions contracts with
time. 
This is an important  point in the program of
rigorously proving the convergence towards the hydrodynamic limit
equation, and  we use this property crucially in our forthcoming work 
\cite{cf:hydro} in the periodic boundary condition setting.
In fact, the idea of the $H^{-1}$ method is to prove that the
$\mathbb L^2$ distance between the deterministic PDE and the randomly
evolving interface stays close to zero at all times. Let us recall
briefly how this works in the Ginzburg-Landau model
\cite{FuSpo,Nishi}. By an entropy production argument
\cite[Th. 4.1]{FuSpo} one shows that at positive times the law of
interface gradients is locally close to a certain equilibrium Gibbs
measure with an unknown slope. The crucial point is that if the slope
is ``wrong'', i.e. different from that of the solution of the PDE, in
which case the random interface has deviated from the deterministic
evolution, the derivative of the $\mathbb L^2$ norm turns out to be
\emph{negative}, which means that the evolution is immediately driven
back to the deterministic one (cf. \cite[Sec. 5.1]{FuSpo}).  In turn, the mathematical mechanism
behind this fact is the same as the one that guarantees that the
$\mathbb L^2$ between two solutions of the PDE contracts with time.


We will also show that the $\mathbb L^1$ distance
between solutions of the limit PDE is non-increasing, and decreases only by a boundary
term (Section \ref{sec:l1}). This is the analogue of the above-mentioned average volume-contraction property of the microscopic dynamics.

Finally let us point out that the exact formula for the hydrodynamic
equation leads to some striking identities involving the surface tension (see notably Eq.
\eqref{eq:muhat} and the discussion in Remark \ref{rem:esotica}) for which it would be very interesting to find a  probabilistic interpretation.

\smallskip

The work is organized as follows. In Section \ref{sec:2} we introduce
precisely the model and the dynamics. The hydrodynamic equation is
given (in two different but equivalent forms) in Section \ref{sec:3}, where
we also discuss some of its properties, notably volume contraction.
In Section \ref{sec:GK} we give a first justification for the limit equation,
based on linear response theory. In Section \ref{sec:LE} instead we derive
the hydrodynamic equation under a local equilibrium
assumption. Finally, in Section \ref{sec:K} we explain how to perform some
useful equilibrium computations.

\section{The model and the dynamics}
\label{sec:2}

\subsection{Monotone surfaces and height function}
\label{sec:monotonesurf}

We start by defining discrete monotone surfaces. 
\begin{Definition}  Let $Q$ be the collection of closed squares
in $\mathbb R^3$ of side $1$, with the four vertices in $\mathbb Z^3$.
  A discrete (or stepped) monotone surface  $\Sigma$ is a
  connected union of elements of $Q$ that projects bijectively on
  the  $P_{111}$ plane (the linear subspace of $\mathbb R^3$ of normal vector $(1,1,1)$).
  
\end{Definition}

\begin{figure}
\includegraphics[width=.5\textwidth]{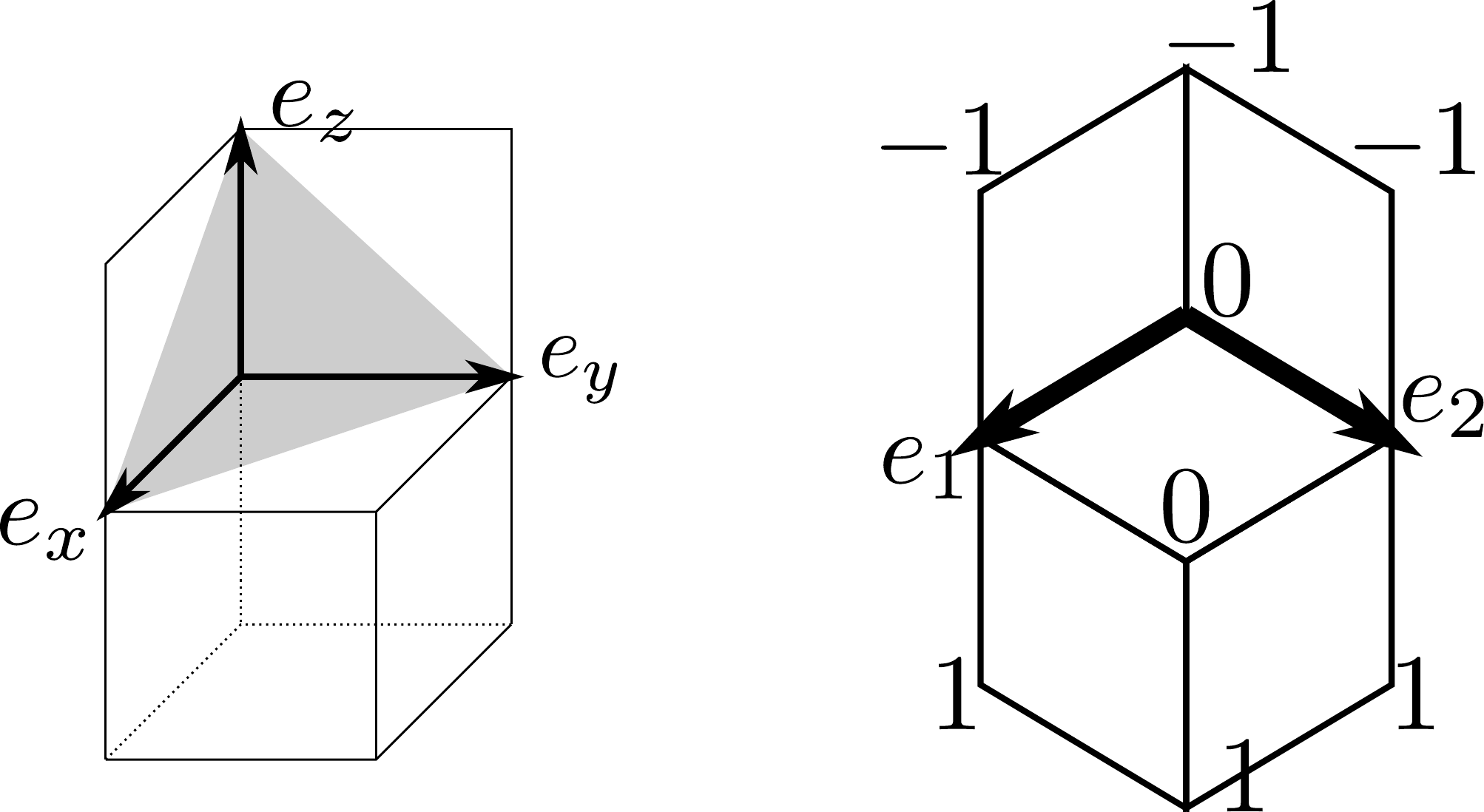}
\caption{Correspondence between monotone surface and lozenge
  tiling. Left: A monotone surface seen as a subset of $\bbR^3$, with
  the axes of $\bbR^3$. The grey triangle lies in
  the plane $P_{111}$. Right: The corresponding lozenge tiling with its
  height function and the vectors $e_1$ and $e_2$.}
\label{fig:B}
\end{figure}

The $P_{111}$ orthogonal projection (denoted $\Pi_{111}$) of a
square face of $\Sigma$ is a lozenge with angles $\pi/3$ and $2\pi/3$,
side-length $\sqrt{2/3}$ and three possible orientations: north-west,
north-east and horizontal, according to whether the normal vector to
the square face in question is $( 1 , 0 , 0 )$, $( 0,1 , 0 )$ or $( 0
, 0 , 1 )$.  The projection of $\Sigma$ gives therefore a lozenge
tiling of $P_{111}$.  Vertices of the lozenges are the vertices of a
triangular lattice $\mathcal T$ of side $\sqrt{2/3}$. We will refer to
north-west oriented, north-east oriented and horizontal lozenges as
lozenges of types $1,2,3$ respectively. See Fig. \ref{fig:B}.

 Let $e_x,e_y,e_z$ be the usual orthonormal vectors of $\mathbb R^3$.
On the plane $P_{111}$ we introduce unit vectors $e_1,e_2$ and
correspondingly coordinates $u=(u_1,u_2)$ as follows: a given
reference vertex $u_0\in\mathcal T$ (for example, the one on which the
origin of $\mathbb Z^3$ projects) has coordinates $(0,0)$ and the
vector $e_1$ (resp. $e_2$), of coordinates $(u_1,u_2)=(1,0)$ (resp.
$(u_1,u_2)=(0,1)$) is the vector from $u_0$ to its nearest neighbor in
direction $e^{-5i \pi/6}$ (resp. $e^{-i \pi/6}$). That is,
$e_1=\Pi_{111}(e_x),e_2=\Pi_{111}(e_y)$. Note that with this choice of
coordinates, triangular faces of $\mathcal T$ have side-length $1$ and not
$\sqrt{2/3}$.

In order to turn the correspondence between stepped surfaces and lozenge
tilings into a bijection, we impose that  $(0,0,0)\in \Sigma$.

Given a discrete monotone surface $\Sigma$, the height function
$h=h_\Sigma:\mathcal T\mapsto \mathbb Z$ is defined as follows:
$h{(u_1,u_2)}$ equals \emph{minus} the height with respect to the
horizontal plane (i.e. \emph{minus} the $z$ coordinate) of the point
$p\in\Sigma$ that projects on $(u_1,u_2)$, i.e. such that
$\Pi_{111}(p)=u_1 e_1+u_2 e_2$. Of course $h(0,0)=0$ since we imposed
$(0,0,0)\in\Sigma$.  The reason for the minus sign is that otherwise
the interface gradients would be given by minus the lozenge densities (see Remark \ref{rem:gradienti} just below),
which would lead to less readable formulas later.  The height function
can be naturally extended to the whole plane by linear interpolation
in each face of $\mathcal T$.

\begin{Remark}
\label{rem:gradienti}
Observe that when one moves by one lattice step in $\mathcal T$ along the $+e_1$ or $+e_2$ directions
the height function increases by $1$ if one
crosses a lozenge, and is unchanged if one moves along the edge of a
lozenge.  When instead one moves by a lattice step upward in the
vertical direction (i.e.  by $-e_1-e_2$), the height function is unchanged if one crosses a
lozenge, and decreases by $1$ if one moves along the edge of a
lozenge.  
\end{Remark}

We will be interested in dynamics in finite domains. For $L=1,2,\dots$
let $U_L$ be a simply connected, bounded union of triangular faces of
$\mathcal T$ that can be tiled by lozenges, and let $\partial U_L$ be
its boundary, seen as a collection of edges of $\mathcal T$. Assume
that the site $(0,0)$ where height is fixed to zero is on $\partial
U_L$. 

Call $\Omega_{U_L}$ the set of lozenge tilings of $U_L$ and $\eta$ the
generic element of $\Omega_{U_L}$. Remark that the height function on
$\partial U_L$ is independent of the configuration $\eta\in
\Omega_{U_L}$.  If wished, one can imagine that $\eta$ is extended to
a lozenge tiling $\eta'$ of the whole plane, just by completing $\eta$
with a tiling $\eta_0$ of $P_{111}\setminus U_L$, fixed once and for
all. From this point of view, $\Omega_{U_L}$ can be identified with
the set of monotone surfaces $\Sigma$ such that $(0,0,0)\in\Sigma$ and
such that the projection of $\Sigma$ restricted to $P_{111}\setminus
U_L$ coincides with $\eta_0$.

We will assume from now on that $U_L$ has a scaling limit in the following sense:
\begin{Assumption}
\label{assUL}
As $L\to\infty$, $(1/L)U_L$ tends in Hausdorff distance to a bounded
simply connected closed domain $U\subset \mathbb R^2$ with smooth
boundary, and the graph of the function
\[
u\in \frac1L \partial U_L\mapsto \frac1L h_{\eta'}(uL)\in \mathbb R
\]
 tends to the graph of a continuous function $\phi:\partial U\mapsto \mathbb R$.
Moreover, there exists a $C^1$  function $\psi:U\mapsto \mathbb R$   such that $\psi=\phi$ on $\partial U$ and that
\begin{eqnarray}
\label{eq:nonmassimale}
  \nabla \psi(u)=(\partial_{u_1} \psi,\partial_{u_2} \psi)(u) \in \mathbb T\quad \forall u\in U
\end{eqnarray}
where 
$\mathbb T\subset \mathbb R^2$ is the open triangle with vertices 
$(0,0),(1,0),(0,1)$.
\end{Assumption}
To understand the condition \eqref{eq:nonmassimale}, recall that
\begin{Theorem}\cite{CKP}
\label{th:CKP}
  Given $\psi:U\mapsto \mathbb R$ satisfying \eqref{eq:nonmassimale} and $\psi|_{\partial U}=\phi$,
  \begin{eqnarray}
    \label{eq:KLP}
   \lim_{\delta\to0}\lim_{L\to\infty}\frac1{L^2}\ln  |\{\eta\in \Omega_{U_L}:\sup_{u\in U}|L^{-1}h_\eta(uL)- \psi(u)|\le \delta\}|= -\int_U \sigma (\nabla \psi) du
  \end{eqnarray}
where $\ln$ denotes the natural logarithm and, for $\rho=(\rho_1,\rho_2)$, 
\begin{eqnarray}
  \label{eq:sigma}
  \sigma(\rho)=\left\{
    \begin{array}{cc}
\frac1\pi\left[
\Lambda(\pi \rho_1)+\Lambda(\pi \rho_2)+\Lambda(\pi(1-\rho_1-\rho_2))\right]\le0,& \rho\in \mathbb T\cup \partial \mathbb T\\
      +\infty & \quad \text{otherwise}
    \end{array}
\right.
\end{eqnarray}
with
\[
\Lambda(\theta)=\int_0^{\theta} \ln (2 \sin(t))dt.
\]
\end{Theorem}
Observe that the argument of $\Lambda$ in \eqref{eq:sigma} is positive, since
$\rho_1,\rho_2,\rho_3:=1-\rho_1-\rho_2\ge0$ if $\rho\in \mathbb T$.
The function $\sigma$ is real analytic and strictly negative in $\mathbb T$,
vanishes when $\rho\in \partial\mathbb T$ and its gradient diverges
when $\partial \mathbb T$ is approached.  The condition
\eqref{eq:nonmassimale} guarantees in particular that the cardinality
of $\Omega_{U_L}$ is exponentially large in $L^2$, i.e. the entropy
per unit area is positive.

In view of the definition of height function, one should think of $\rho_i$ as the density of lozenges of types $i=1,2,3$. We however emphasize that making point-wise sense of this intuition is a delicate problem.

\subsection{Translation-invariant Gibbs states}
\begin{figure}[h]
  \includegraphics[width=10cm]{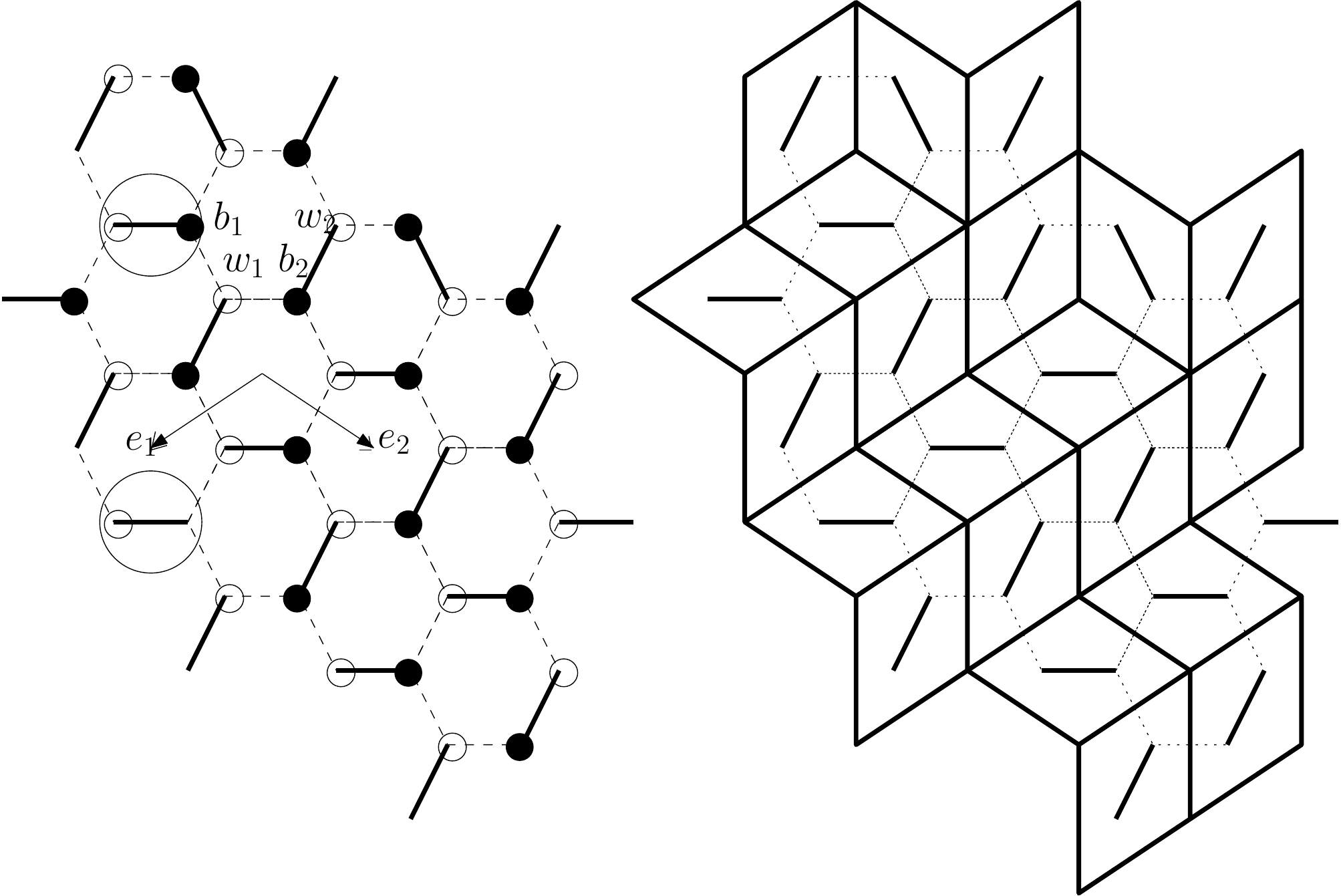}
\caption{The bijection between perfect matching of $\mathcal H$ (left) and lozenge tiling (right).
The vertices $w_1,b_1,w_2,b_2$ have coordinates respectively $(0,0),(0,-1),(-1,0),(0,0)$.
The two encircled horizontal dimers (``particles'') have the same column coordinate $u_1-u_2$ and vertical coordinate differing by $2$.}
\label{fig:1}
\end{figure}

It is well known \cite{Kenyonnotes} that for every $\rho\in \mathbb T$
there exists an unique translation invariant, ergodic Gibbs state on
the set of lozenge tilings of the plane, such that the density of
lozenges of type $i$ is $\rho_i$.  Such measures have the following
explicit form. Let the hexagonal lattice $\mathcal H$ be the dual of
$\mathcal T$ and color its vertices black/white in an alternate way,
see Figure \ref{fig:1}.  A white vertex is given the same coordinate
$x=(x_1,x_2)$ as the black vertex just to its right, and the
coordinates $x_1,x_2$ correspond to the axes $e_1,e_2$ introduced
above. There exists a natural bijection between lozenge tilings $\eta$
of the plane and perfect matchings $M$ of $\mathcal H$, see Figure \ref{fig:1}.

Take a triangle with angles $\theta_i=\pi \rho_i, i=1,2,3$ and let
$k_i$ be the length of the side opposite to $\theta_i$. Given an edge
$e$ of $\mathcal H$, say that it is of type $1,2$ or $3$ if it is
north-west, north-east oriented or horizontal and let $K(e):=k_i$.
Then, given an integer $n$ and edges $e_i,i\le n$ such that the white
(resp. black) vertex of $e_i$ is $w_i$ (resp. $b_i$), one has
\cite{Kenyonnotes}
\begin{eqnarray}
\label{eq:piro}
  \pi_\rho(e_1, \ldots, e_n \in M)=K(e_1)\dots K(e_n)\det(\{K^{-1}(w_{i},b_j)\}_{i,j=1,\dots,n})
\end{eqnarray}
where, if $w$ has coordinates $(x_1,x_2)$ and $b$ has coordinates $(y_1,y_2)$,
\begin{eqnarray}
  \label{eq:K-1}
  K^{-1}(w,b)=\frac1{(2\pi i)^2}\oint \frac{dz}z \oint \frac{dw}w  \frac{z^{x_2-y_2}w^{y_1-x_1}}{k_3+k_1 z+k_2 w}
\end{eqnarray}
and the integrals runs over the torus $\{(z,w)\in \mathbb C^2:
|z|=|w|=1\}$.  Note that \eqref{eq:piro} is unchanged if all $k_i$ are
multiplied by a common factor.  In particular, if $e=(w,b)$ is an edge
of type $i$ one has
\begin{eqnarray}
\label{eq:K1i}
  \rho_i=k_i K^{-1}(w,b).
\end{eqnarray}

\subsection{Particles} We have seen that the $P_{111}$ projection of any discrete monotone surface gives a lozenge tiling $\eta$ of the plane. 
 Horizontal lozenges will be  called ``particles'' and will be given a
 label $b$. To each particle will be associated a ``vertical position''
 $n(b)$, defined as
\begin{eqnarray}
  \label{eq:nb}
  n(b)=\frac{u_1+u_2}2,
\end{eqnarray}
with $(u_1,u_2)$ the coordinates of the upper corner of the particle
(horizontal lozenge),
as well as a ``horizontal position'' (or ``column coordinate'') 
\[
c(b)=u_1-u_2,
\] 
see Figure \ref{fig:1}.
Note that $n(b)\in \mathbb Z$ if $c(b)\in 2\mathbb N$
(i.e. if the column containing $b$ has the same parity as the column containing the
vertex $(0,0)$), and $n(b)\in \mathbb Z+1/2$ otherwise.

Recalling Remark \ref{rem:gradienti}, observe that when the vertical
coordinate of a particle $b$ changes by $\pm |n|$, there are $n$ vertices
in the triangular lattice $\mathcal T$ where the height changes by
$\pm1$.

\medskip

It is well known (and easy to check) that a lozenge tiling of the
plane is uniquely determined by the particle positions, provided that
there is at least one particle per column, which we will assume
henceforth. Recall that the height function on $\partial U_L$ is independent of the configuration $\eta\in \Omega_{U_L}$. From the definition of height function, we deduce that for each column $i$, the number of particles on column $i$ that are in $U_L$ is the same for every  $\eta\in \Omega_{U_L}$. Actually the whole tiling $\eta$ is
uniquely determined (once $U_L$ and $\eta_0$ on $P_{111}\setminus U_L$
are given) by the positions of the particles in $U_L$.

It is also well known and easy to check that particle positions satisfy the
following interlacement properties: if $b,b'$ are two particles on the
same column $i=c(b)=c(b')$ with $n(b)<n(b')$ and if there is no particle $b''$ on
column of index $i$ with $n(b)<n(b'')<n(b')$ then there is
a unique particle $b^{right}$ (resp. $b^{left}$) in column $i+1$
(resp. $i-1$) such that $n(b)<n(b^{right})<n(b')$ (resp. $n(b)<n(b^{left})<n(b')$).

\begin{figure}[h]
  \includegraphics[width=11cm]{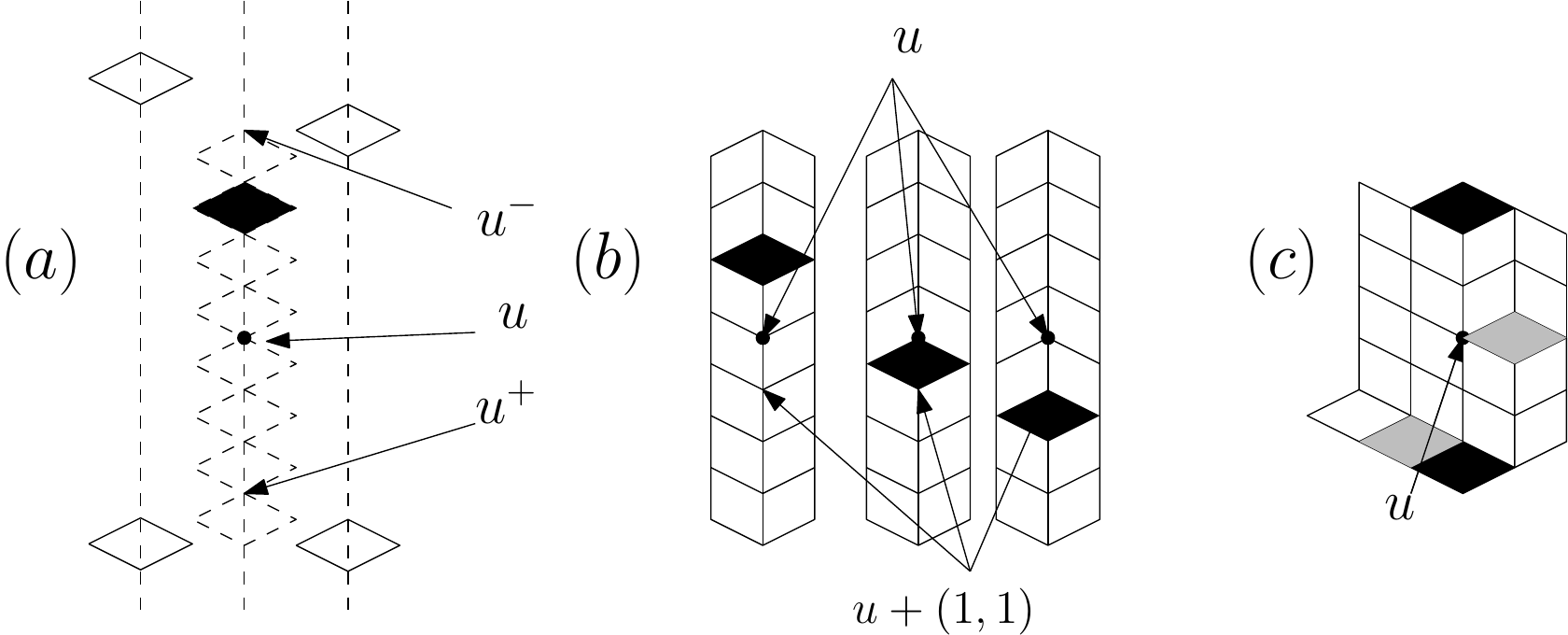}
  \caption{(a): For the black particle $b$ one has
    $n^+(b)=(u^+_1+u^+_2)/2$ and $n^-(b)=(u^-_1+u^-_2)/2$ (recall that
    vertical positions are measured w.r.t. $-e_1-e_2$, which is why
    $u^+$ appears below $u^-$ in the drawing). Further moves are
    prevented by the four white particles. The event $X(u)$ is
    verified here and particle $b^{(u)}$ is the black one. (b):
    According to the position of $b^{(u)}$ w.r.t. $u$, the edge
    $(u,u+(1,1))$ can either cross the particle $b^{(u)}$ (mid
    drawing) or be the common edge of two lozenges, of types $1$ and
    $2$ (left and right drawing). (c) In this configuration, the event
    $X(u)$ is not realized. None of the particles (drawn in black) in
    the column of $u$ can take vertical position $(u_1+u_2)/2$ without
    pushing particles (in gray) in neighboring columns.}
\label{fig:Xu}
\end{figure}

In the study of our interface dynamics we will need the following two definitions:
\begin{Definition}
  Given $\eta\in \Omega_{U_L}$ and a particle $b$ in $U_L$, we let
  $n^+(b), n^-(b)$ be the largest/smallest possible vertical position
  that particle $b$ can take in any configuration $\eta'\in\Omega_L$ such that
  all particles other than $b$ have the same position as in $\eta$. See Figure \ref{fig:Xu} (a).
Of course $n^\pm(b)$ is a function of $\eta$.  Also, we call
  $I(b)=\{n^-(b),\dots,n^+(b)\}$. 
\end{Definition}

\begin{Definition}
\label{def:Xu}
Given a lozenge tiling of the plane and
$u=(u_1,u_2)\in\mathcal T$
we
call $X(u)$ the event that there exists a particle 
$b^{(u)}$ with column coordinate
$u_1-u_2$ such that \[n^-(b^{(u)})\le \frac{u_1+u_2}2\le n^+(b^{(u)}), \] 
see Figure \ref{fig:Xu} (a). (Such particle is necessarily unique). The particle $b^{(u)}$ can be moved to position $(u_1+u_2)/2$ without moving any other particle.

The event $X(u)$ is equivalently
the event that the edge $u,u+(1,1)$ of $\mathcal T$ is either the common side of two lozenges, one of type $1$ and one of type $2$, 
or it crosses  a lozenge of type $3$, see Figure \ref{fig:Xu} (b). In the latter case, one has $n(b^{(u)})=(u_1+u_2)/2$.
Therefore, the event $X(u)$ is not realized when the edge  $u,u+(1,1)$ is the common side of two lozenges, both of type $1$ or both of type $2$ (e.g. Fig. \ref{fig:Xu} (c)).

\end{Definition}

\subsubsection{Interpretation of $\rho$ in terms of particles}

One can give an interpretation, purely in terms of interlaced
particles, to the two parameters $\rho_1,\rho_2$ labelling the Gibbs
measures $\pi_\rho$. First of all, $1-\rho_1-\rho_2$ is the density of
particles in any given column. The difference $\rho_1-\rho_2$
corresponds to an asymmetry parameter as follows. Look at a column, say the one labelled $0$, and call $\{n_i\}_{i\in\mathbb Z}$ the vertical
positions of its particles $\{b_i\}_{i\in\mathbb Z}$, ordered so that $n_i<n_{i+1}$. Given
particles $b_i$ and $b_{i+1}$, let $b'_i$ be the unique particle in column
$1$ whose vertical position $n'_i$ satisfies $n_i<n'_i<n_{i+1}$.  Then
one has
\begin{eqnarray}
  \label{eq:18}
 \frac{ \rho_1-\rho_2}2=\lim_{K\to\infty}\frac{\sum_{i=1}^K n'_{i}-\frac{n_{i+1}+n_i}2}{\sum_{i=1}^K(n_{i+1}-n_i)}
\end{eqnarray}
where the limit holds $\pi_\rho$-almost surely, due to ergodicity of
the Gibbs measure. In other words, $\rho_1-\rho_2$ is a measure of how
much $n'_i$ is biased away from the mid-point $(n_i+n_{i+1})/2$.

\begin{figure}[h]
  \includegraphics[width=2.5cm]{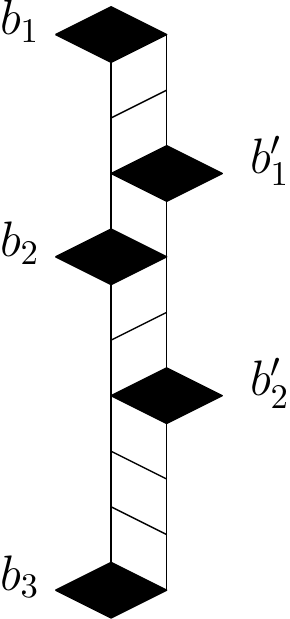}
\caption{If for instance particle $b_1$ has position $n_1=0$, then
  $n_2=4$, $n_3=10$ and $n'_1=5/2$, $n'_2=13/2$. The number of type-1
  lozenges is $N^{(2)}_1=4=(n_2-n'_1-1/2)+(n_3-n'_2-1/2)$.}
\label{fig:bipi}
\end{figure}

To see why \eqref{eq:18} holds, look at Figure \ref{fig:bipi}: running along
column $0$ from position $n_1$ to $n_{K+1}$, 
 the number of lozenges of
type $1$ that are adjacent to column $0$ to its right is 
\begin{eqnarray}
  \label{eq:19}
 N^{(K)}_1= \sum_{i=1}^K (n_{i+1}-n'_i-1/2)
\end{eqnarray}
and the number of lozenges of type $2$ is 
\begin{eqnarray}
  \label{eq:20}
 N^{(K)}_2= \sum_{i=1}^K (n'_i -n_{i}-1/2).
\end{eqnarray}
The factors $1/2$ keep into account the fact that particle positions
in column $0$ are integers and those in column $1$ are half-integers.
One has then 
\begin{eqnarray}
  \label{eq:21}
  \frac{\rho_1-\rho_2}{\rho_1+\rho_2}=\lim_K
  \frac{N^{(K)}_1-N^{(K)}_2}{N^{(K)}_1+N^{(K)}_2}=
\lim_K\frac{\sum_{i=1}^K(n_{i+1}+n_i-2n'_i)}{\sum_{i=1}^K(n_{i+1}-n_i-1)}.
\end{eqnarray}
On the other hand, 
\begin{eqnarray}
  \label{eq:22}
  \frac{\sum_{i=1}^K(n_{i+1}-n_i-1)}{\sum_{i=1}^K(n_{i+1}-n_i)}=1-\frac
  K {n_{K+1}-n_1}
\end{eqnarray}
that converges to $1-\rho_3=\rho_1+\rho_2$ since $K$ is
the number of particles in column $0$ in a segment of length
$n_{K+1}-n_1$. Equation \eqref{eq:18} then follows.

\subsection{Dynamics}

We will study two continuous-time Markov dynamics on
$\Omega_{U_L}$. Both are reversible with respect to the uniform
measure $\mathbb P_{U_L}$. We will only define the dynamics in terms of movements of particles but recall that these determine the whole tiling. In the dynamics, only particles in $U_L$
can evolve.

\begin{Definition}[Dynamics I]
\label{def1}
  For any particle  $b$ in $U_L$ and any   $k\in I(b),k\ne n(b)$,
$b$ moves to vertical  position $k$ with rate $1/(2|n(b)-k|)$.
\end{Definition}
This is equivalent to a dynamics introduced by Luby, Randall and Sinclair \cite{LRS}. Let us recall that this dynamics can be used  
 as an auxiliary process to show that  the ``single-flip'' dynamics,
 where particles are instead allowed to move only to $n(b)\pm1$ with
 equal rates (provided $n(b)\pm 1\in I(b)$), has a mixing time and inverse spectral gap that is at most polynomial in $L$.

\begin{Definition}[Dynamics II]
\label{def2}
  For any particle $b$ and any position  $k\in I(b),k\ne n(b)$,
$b$ moves to position $n$ with rate $1/|I(b)|$.
\end{Definition}
In other words, with rate $1$ each particle is redistributed uniformly
among its instantaneously available positions. This dynamics was
introduced in \cite{CMST}, again as an auxiliary process to analyze
the single-flip dynamics.

It is immediate to see that both are reversible w.r.t the uniform measure. We will call $\mathcal L^{I}, \mathcal L^{II}$ the generators of the two dynamics.
The configuration at time $t$ will be denoted $\eta(t)$ and dependence
on the boundary condition as well as the index $I,II$, that distinguishes 
between the two dynamics will not be indicated explicitly.

\begin{Remark}
\label{rem:asymm}
In \cite{Toninelli2+1} was defined an irreversible (driven) dynamics
on lozenge tilings of the infinite triangular lattice $\mathcal T$. In
its totally asymmetric version, each particle $b$ jumps to any $k\in
I(b), k>n(b)$ with rate $1$. It was proven in \cite{Toninelli2+1} that
the Gibbs measures $\pi_\rho$ is stationary for such driven dynamics
(a consistent part of the work consisting in proving that the dynamics
is well defined for almost every initial condition sampled from
$\pi_\rho$).  We will see (cf. discussion just after Theorem \ref{prop:K}) that the two reversible dynamics of
Definitions \ref{def1} and \ref{def2} are not unrelated to the
irreversible dynamics of \cite{Toninelli2+1}.

\end{Remark}

\section{The limit hydrodynamic equation}
\label{sec:3}

Call $\eta=\eta(0)$ the initial condition of the dynamics (actually $\eta$ is a sequence $\{\eta_L\}_{L\in \mathbb
N}$, but we drop the subscript $L$) and $\eta(t)$ the
configuration at time $t$ (recall that we do not distinguish between dynamics I
and II in the notation).   Assume that $\eta$ approximates a smooth profile, i.e. there exists $\psi_0$ satisfying \eqref{eq:nonmassimale} such that
\begin{eqnarray}
  \label{eq:iniziale}
  \lim_{L\to\infty}\frac1L h_{\eta}( uL)=\psi_0(u)
\end{eqnarray}
for every $u\in U$. 
Let for $t\ge0, u\in U$
\[
 H(u,t)=\frac1L h_{\eta(t L^2)}( uL).
\]
On general grounds \cite{Spohn}, one expects $H$ to concentrate around
a deterministic solution in the sense that there exists
some deterministic function $\{\psi(u,t)\}_{u\in U, t\ge0}$  such that, for every $\epsilon>0$ and $t\ge0$,
\begin{eqnarray}
  \label{eq:generale}
  \mathbb P(\exists u\in U: |H(u,t)-\psi(u,t)|\ge \epsilon)\stackrel{L\to\infty}\to 0.
\end{eqnarray}
 Furthermore, $\psi$ should follow a  non-linear PDE of the form
 \begin{eqnarray}
   \label{eq:PDE}
 \left\{
   \begin{array}{ll}
   \partial_t \psi=\mu(\nabla \psi)\sum_{i,j=1,2}\sigma_{i,j}(\nabla \psi)\frac{\partial^2}{\partial_{u_i}\partial_{u_j}}\psi\\   
  \psi(u,t)=\psi_0(u) & \text{ if } t=0 \text{ or if } u\in \partial U
   \end{array}
 \right.
 \end{eqnarray}
 where $\sigma_{i,j}(\rho):=\partial^2\sigma/\partial \rho_i\partial
 \rho_j$, with $\sigma$ defined in \eqref{eq:sigma}, and $\mu(\cdot)$
 is a positive function. This equation is of parabolic type, since the
 Hessian matrix $\big(\sigma_{i,j}(\rho)\big)_{i,j=1,2}$ is strictly
 positive definite for $\rho\in\mathbb T$ (positive definiteness
 follows from convexity of the surface tension $\sigma$ and strict
 positivity follows from the fact that the determinant of
 $\big(\sigma_{i,j}(\rho)\big)_{i,j=1,2}$ equals identically $\pi^2$
 \cite[Th. 5.5]{KOS}, as can also be checked from \eqref{eq:3} below).

 The positive coefficient $\mu(\cdot)$ is called the ``mobility'' and in
 general will depend on the microscopic definition of the dynamics. In
 particular, a priori there is no reason for it to be the same for
 dynamics I and II, but we will see below that the mobility does in
 fact coincide in the two cases. 
 \begin{Remark}
   While we see no a priori physical reason why the two mobilities
   should coincide, let us comment on what is behind this equality.
   Both dynamics satisfy a form of ``gradient condition'', that is
   responsible for the vanishing, thanks to an exact summation by
   parts, of the term in the Green-Kubo expression \eqref{eq:muu} that
   involves equilibrium correlations at different times. Then,
   $\mu(\rho)$ turns out to be proportional to the derivative at time
   zero of the mean square displacement of a particle, for the process
   started from the stationary state $\pi_\rho$:
\begin{eqnarray}
  \label{eq:smi}
\mu(\rho)\propto \pi_\rho\left(\sum_{y:n(b)+y\in I(b)} c_{b,y}(\eta)y^2\right),  
\end{eqnarray}
with $c_{b,y}(\eta)$ the rate at which particle $b$ jumps from
position $n(b)$ to $n(b)+y$ in configuration $\eta$.  See the first
term in \eqref{eq:muu}. A simple and explicit computation, using the
fact that at equilibrium the position of a particle is uniform given
the position of the other particles, shows that, despite the fact that
the transition rates are different for dynamics I and II, the average
in \eqref{eq:smi} is the same in both cases (see
Eq. \eqref{eq:summa}).  It would be interesting to understand whether
there are other natural transition rates that lead to the same
mobility.

Another hint that the two dynamics have common features is that both were designed (in \cite{LRS} and \cite{CMST} respectively) to have the property of contracting the mutual volume between two configurations. 
 \end{Remark}
 
The meaning of \eqref{eq:PDE} is that, in the diffusive
 scaling, the interface velocity will be given by the gradient flow
 associated to the surface tension functional
\begin{eqnarray}
  \label{eq:6}
  F(\psi)= \int_U  \sigma(\nabla \psi) du
\end{eqnarray}
times a certain mobility
coefficient $\mu$ that depends on the local slope. Note indeed that
\begin{eqnarray}
  \label{eq:7}
\sum_{i,j=1,2}\sigma_{i,j}(\nabla \psi)\frac{\partial^2}{\partial_{u_i}\partial_{u_j}}\psi=-
\frac{\delta F(\psi)}{\delta \psi(u)}.  
\end{eqnarray}
More explicitly, one finds from 
\eqref{eq:sigma}
\begin{equation}
  \label{eq:3}
\begin{aligned}
\sigma_{i,i}(\rho)&=
\pi\cot(\pi \rho_i)+\pi\cot(\pi(1-\rho_1-\rho_2))
\\
\sigma_{1,2}(\rho)&=\pi\cot(\pi(1-\rho_1-\rho_2)).
\end{aligned}
\end{equation}

An expression for $\mu(\rho)$ can be obtained from linear response
theory.  Usually (cf. \cite[Section 4]{Spohn} and \eqref{eq:muu}
below), such expression is given by the sum of two terms: the first
involves the average w.r.t. $\pi_{\rho}$ of a local observable and the
second involves the integral over time $t$, ranging from $0$ to
$\infty$, of the correlations (in the stationary process started from
$\pi_{\rho}$) between an observable at time $0$ and an observable at
time $t$.  In general, it is not possible to compute the second term
explicitly. In lucky  cases (e.g. the zero-temperature dynamics of
interfaces of the 2D Ising model or the Langevin dynamics of the
Ginzburg-Landau effective interface model \cite{Spohn}) the second
term vanishes due to a summation by parts. This turns out to be the case also for our dynamics.

The expression for the mobility, provided by linear response theory,
is in our case
(see Section \ref{sec:GK} for a  derivation along the lines of \cite{Spohn})
\begin{eqnarray}
  \label{eq:mu}
  \mu(\rho)=\frac1{12}\pi_{\rho}\left[ (|I(b^{(0,0)})|+1)1_{X(0,0)}1_{n(b^{(0,0)})\ne 0}\right],
\end{eqnarray}
for both dynamics I and II.  Here, $b^{(0,0)}$ and $X(0,0)$ are just
$b^{(u)},X(u)$ as in Definition \ref{def:Xu}, and we
have arbitrarily chosen $u=(0,0)$ by translation invariance.

In Section \ref{sec:K} we show:
\begin{Theorem}
\label{prop:K}
 The r.h.s. of \eqref{eq:mu} equals
\begin{eqnarray}
  \label{eq:CF}
\frac{V(\rho)}2:=\frac12\pi_\rho\left[(n^+(b^{(0,0)})+1)1_{X(0,0)}1_{n(b^{(0,0)})<0}\right].
\end{eqnarray}
  
\end{Theorem}
The advantage of the rewriting \eqref{eq:CF} is that $V(\rho)$ turns out to be nothing but the average interface velocity for
the totally asymmetric process defined in Remark \ref{rem:asymm}, in the stationary
measure $\pi_\rho$. Namely, let the ``total current'' $J(0,t)$ denote the
number of particles that cross a fixed vertex of the triangular lattice
$\mathcal T$ in the time interval $[0,t]$ for the asymmetric process.  Then \cite{Toninelli2+1}
\begin{eqnarray}
  \label{eq:2}
  \langle J(0,t)\rangle_{\rho}=t V(\rho),
\end{eqnarray}
where $\langle\cdot\rangle_{\rho}$ denotes average w.r.t. the
stationary process started from $\pi_\rho$.
In \cite[Th. 2.7]{FerrariSunil} it was proven (with somewhat different notations)
that 
\begin{eqnarray}
  \label{eq:1}
{V(\rho)}=\frac1{\pi} \frac{\sin(\pi \rho_1)\sin(\pi \rho_2)}{\sin(\pi(1-\rho_1-\rho_2))}
\end{eqnarray}
Recall that $0< \rho_1,\rho_2,(1-\rho_1-\rho_2)<1$ and that these three numbers
give the average fraction of 
lozenges  of types $1,2,3$ respectively under the measure $\pi_{\rho}$.
In conclusion, both for dynamics I and II, the 
linear response theory mobility defined as in
\eqref{eq:mu} equals 
\begin{eqnarray}
  \label{eq:4}
  \mu(\nabla \psi)=\frac1{2\pi} \frac{\sin(\pi \partial_{u_1}\psi)\sin(\pi \partial_{u_2}\psi)}{\sin(\pi(1-\partial_{u_1}\psi-\partial_{v_2}\psi))}.
\end{eqnarray}
The conjectural explicit form of the hydrodynamic limit equation is then given by   \eqref{eq:PDE}, together with \eqref{eq:3} and \eqref{eq:1}.

\subsection{Hydrodynamic equation and volume contraction}
\label{sec:vd}
It goes beyond the scopes of the present work to investigate the existence  and 
regularity of the solutions of \eqref{eq:PDE}. This might be a
non-trivial issue due to the singularity of $\mu(\cdot)$ and
$\sigma_{i,j}(\cdot)$ when their argument approaches $\partial\mathbb
T$.  In the following of this section, we will implicitly assume that
the domain $U$ and the initial condition $\psi_0(\cdot)$ are regular
enough that \eqref{eq:PDE} admits a unique classical solution $\psi(u,t)$
that is $C^1$ in $ U\times [0,\infty) $ where we recall that the
domain $U$ is closed. In the forthcoming \cite{cf:hydro} we explain how to extract such existence, uniqueness and smoothness statements from the existing literature (e.g. \cite[Chap. XII]{Liberman}).

It is interesting to remark that \eqref{eq:PDE} can be rewritten as
follows:
\begin{eqnarray}
\label{eq:PDEdiv}
 \partial_t \psi={\rm div} (W\circ \nabla \psi):=\partial_{u_1}W^{(1)}(\partial_{u_1}\psi,\partial_{u_2}\psi)+\partial_{u_2}W^{(2)}(\partial_{u_1}\psi,\partial_{u_2}\psi)
\end{eqnarray}
where $W=(W^{(1)},W^{(2)})$ and
\begin{eqnarray}
  \begin{aligned}
  W^{(1)}(\rho_1,\rho_2)&=-\frac1{2\pi}\cot(\pi(\rho_1+\rho_2))\sin^2(\pi\rho_2)-\frac{\rho_2}4+\frac1{4\pi}\sin(2\pi\rho_2)\\
W^{(2)}(\rho_1,\rho_2)&=W^{(1)}(\rho_2,\rho_1).    
  \end{aligned}
\end{eqnarray}
(This can be checked via a direct computation, using the definition of $W,$ of $\mu$ and
the expressions \eqref{eq:3} for $\sigma_{i,j}$).   One can also
  check that the curl of the vector field $W$ is non-zero, which prevents from writing $W$ as the gradient of 
some function $\hat \sigma$ and the equation \eqref{eq:PDEdiv} for $\psi$ as the gradient flow w.r.t. the associated functional 
$-\int_U \hat \sigma(\nabla \psi)du$.

The rewriting \eqref{eq:PDEdiv} has two interesting consequences, namely contractions in time of
both the $\mathbb L^1$ and $\mathbb L^2$ distances between solutions. The two phenomena are
somewhat different: as we see in a moment, $\mathbb L^1$ contraction is only a
boundary effect, while $\mathbb L^2$ contraction is a bulk effect.

\subsubsection{$\mathbb L^1$ contraction}
\label{sec:l1}

By Gauss' theorem, \eqref{eq:PDEdiv} implies that the
time derivative of the total volume below the surface, $\mathcal V_t(\psi(\cdot,t)):=\int_U
\psi(u,t)du$, is only a boundary term:
\begin{eqnarray}
  \label{eq:dvol}
  \frac d{dt}\mathcal V(\psi(\cdot,t))= \int_{\partial U} W(\nabla\psi(u,t))\cdot n\; d\gamma,
\end{eqnarray}
with $n$ the exterior normal vector to $\partial U$. A stronger property holds. Let  
$\psi^{(1)},\psi^{(2)} $ be two smooth initial conditions for
\eqref{eq:PDEdiv} with \[\psi^{(1)}(u)\ge \psi^{(2)}(u) \text{\; for every\;}
u\in U\quad\text{and} \quad
\psi^{(1)}|_{\partial U}=\psi^{(2)}|_{\partial U}=\psi_0|_{\partial U}.
\] Then, one can show (see end of Section \ref{sec:st2}) that 
\begin{eqnarray}
\label{eq:driftnegativo}
  \int_{\partial U} [W(\nabla\psi^{(1)})-W(\nabla\psi^{(2)})]\cdot n\; d\gamma\le0.
\end{eqnarray}
 Inequality $\psi^{(1)}(u,t)\ge \psi^{(2)}(u,t)$ remains true for
all times, by the usual comparison principle for parabolic PDEs
\cite[Ch. 3]{Protter} (another way to convince oneself that order is preserved is to recall
 that the microscopic dynamics is monotone \cite{Wilson,CMT}: if two initial conditions $\eta,\eta'$
satisfy $h_\eta\le h_{\eta'}$ everywhere, then the two evolutions can be coupled in a way that domination is preserved at all times).
One
concludes that
 \begin{eqnarray}
   \frac d{dt}(\mathcal V(\psi^{(1)}(\cdot,t))-\mathcal V(\psi^{(2)}(\cdot,t)))  =\int_{\partial U} [W(\nabla\psi^{(1)}(u,t))-W(\nabla\psi^{(2)}(u,t))]\cdot n\; d\gamma\le0:
 \end{eqnarray}
 the drift of mutual volume is a boundary effect and is negative. This fact has a
microscopic analog: in fact, if $\eta^{(1)},\eta^{(2)}$ are two configurations in $\Omega_{U_L}$ with 
$h_{\eta^{(1)}}\ge h_{\eta^{(2)}}$ everywhere in $U_L$, then 
the volume drift
\begin{eqnarray}
  \label{eq:Ldri}
  [\mathcal L\sum_{u\in U_L}h(u)](\eta^{(1)})-  [\mathcal L\sum_{u\in U_L}h(u)](\eta^{(2)})
\end{eqnarray}
(with $\mathcal L$ the generator of the Markov chain) is negative and
is non-zero only due to a boundary effect.  This was proven in
\cite{LRS} for dynamics I and in \cite{CMST} for dynamics II, and is
actually the crucial step in the proof that the mixing time is
polynomial in $L$. In this perspective, it is natural to recover such
volume decrease property in the hydrodynamic equation.

It is worth emphasizing that volume contraction is not an a-priori obvious property. In particular, it is easy to check that the single-flip dynamics, at the microscopic level,  \emph{does not} contract volume.

\subsubsection{$\mathbb L^2$ contraction} \label{sec:l2} Let again $\psi^{(1)}(u,t),
\psi^{(2)}(u,t)$ be two smooth solutions of \eqref{eq:PDE}, with the same
boundary data on $\partial U$. This time we do not require that $\psi^{(1)}\ge \psi^{(2)}$. We have
\begin{multline}
  \label{eq:13}
  \frac d{dt}\int_U (\psi^{(1)}(u,t)-\psi^{(2)}(u,t))^2du\\=-2\int_U(\nabla \psi^{(1)}(u,t)-\nabla \psi^{(2)}(u,t))\cdot(W(\nabla \psi^{(1)}(u,t))-W(\nabla \psi^{(2)}(u,t)))du.
\end{multline}
We claim now that
\begin{eqnarray}
  \label{eq:23}
  (a-b)\cdot(W(a)-W(b))\ge0
\end{eqnarray}
whenever $a,b$ belong to the triangle $\mathbb T$, which implies that
the time derivative in \eqref{eq:13} is negative.
To prove \eqref{eq:23}, it is sufficient to prove that the matrix
\begin{eqnarray}
  \label{eq:24}
  H_W(\rho):=\left(
  \begin{array}{cc}
    \partial_{\rho_1} W^{(1)}(\rho) & \partial_{\rho_2} W^{(1)}(\rho)\\
\partial_{\rho_1} W^{(2)}(\rho) & \partial_{\rho_2} W^{(2)}(\rho)
  \end{array}
\right)
\end{eqnarray}
is positive definite for every $\rho\in \mathbb T$.
The trace of $H_W$ is
\begin{eqnarray}
  \label{eq:25}
  {\rm Tr}(H_W(\rho))=\frac12 \frac{\sin(\pi \rho_1)^2+\sin(\pi \rho_2)^2}{\sin(\pi(\rho_1+\rho_2))^2} >0
\end{eqnarray}
(recall that $\rho_1,\rho_2\in (0,1)$)
while
\begin{multline}
  \label{eq:26}
  \det (H_W(\rho))= \frac1{64\sin(\pi(\rho_1+\rho_2))^3} \left\{
5 \sin[\pi (\rho_1 + \rho_2)] + 
\sin[3\pi (\rho_1 + \rho_2)]\right.\\\left.- 
   2 \sin[\pi (3\rho_1 + \rho_2)]-2\sin[\pi (\rho_1 + 3\rho_2)]\right\}.
\end{multline}
Given that $0<\rho_1+\rho_2<1$ for $\rho\in\mathbb T$, the ratio in
the r.h.s. of \eqref{eq:26} is positive. As for the term $\{\dots\}$, 
one can check that its unique extremum for
$\rho\in \mathbb T$ is at $\rho=(\arctan(\sqrt5)/\pi,\arctan(\sqrt 5)/\pi)$, where $\{\dots\}>0$.

\begin{Remark}
If the interface mobility $\mu(\cdot)$ were a constant
(i.e. if $W$ were proportional to $\nabla \sigma$, in which case 
\eqref{eq:PDE} would  be the gradient flow w.r.t. the surface tension functional) then $\mathbb L^2$
contraction would be an immediate consequence of convexity of the
surface tension, since the matrix $H_W(\rho)$ would be replaced by the
Hessian of $\sigma(\sigma)$ computed at $\rho$.
\end{Remark}

\subsection{Another form for the hydrodynamic equation}

There exists another way of guessing the
hydrodynamic limit equation, this time not based on linear response
but on a ``local equilibrium'' assumption. For this derivation, it is
actually more convenient to use a different way of parametrizing the
interface and the height function. 
The corresponding expression for the hydrodynamic limit equation will show an interesting
link between mobility and surface tension, see \eqref{eq:muhat}.

\subsubsection{Level set function}

Let $P_{110}$ be the linear subspace of $\mathbb R^3$
 orthogonal to $(1,1,0)$ (i.e. the plane $x+y=0$). On $P_{110}$ we take
coordinates $v=(v_1,v_2)$ whose unit vectors $\hat e_1,\hat e_2$ are the
$P_{110}$ orthogonal projections of the Cartesian unit vectors
$e_y,e_z$ of $\mathbb R^3$ (i.e. $\hat e_1=\Pi_{110}(e_y), \hat
e_2=\Pi_{110}(e_z)$) and such that the point of coordinates $v=(0,0)$ is the $P_{110}$ projection of $(0,0,0)$.
Given $v=(v_1,v_2)\in \mathbb Z\times (\mathbb Z+1/2)$, let 
\begin{eqnarray}
  \label{eq:hatac}
  \hat h(v):= -p \cdot (e_x+e_y),
\end{eqnarray}
with $ \cdot $ the usual scalar product on $\mathbb R^3$
and $p\in \mathbb Z^3$ the unique\footnote{The reason for the
  choice $v_2\in\mathbb Z+1/2$ is that if instead we required that
  $v_2\in \mathbb Z$, then $p$ would not be uniquely defined
  (horizontal square faces of $\Sigma$ are projected into segments at some
  integer vertical height $v_2$).} point of the surface $\Sigma$ whose
$\Pi_{110}$ projection is $(v_1,v_2)$.
We can easily extend $\hat h$ to a function on $\mathbb R^2$: formula
\eqref{eq:hatac} is well defined whenever $v_2\not\in \mathbb Z$, and
for $v_2\in \mathbb Z$ (in which case the point $p$ may be not
uniquely defined) we let, say,  $\hat h(v)=\lim_{\epsilon\searrow 0} \hat
h(v_1,v_2-\epsilon)$. This choice is somewhat arbitrary but this should be irrelevant in the $L\to\infty$ limit.
As in the case of the height function $\{h(u)\}_{u\in \mathcal T}$, the function
$\{\hat h(v)\}_{v\in \mathbb Z\times (\mathbb Z+1/2)}$ uniquely
determines the surface $\Sigma$.

\begin{Remark}
\label{rem:duevolumi}
Note that when a particle moves one step up/down, the height
function $h$ changes by $-1/+1$ at some vertex $u\in \mathcal T$, while the function $\hat h$ changes by $-2/+2$ at some $v\in
\mathbb Z\times (\mathbb Z+1/2)$.
\end{Remark}

\begin{figure}[h]
  \includegraphics[width=6cm]{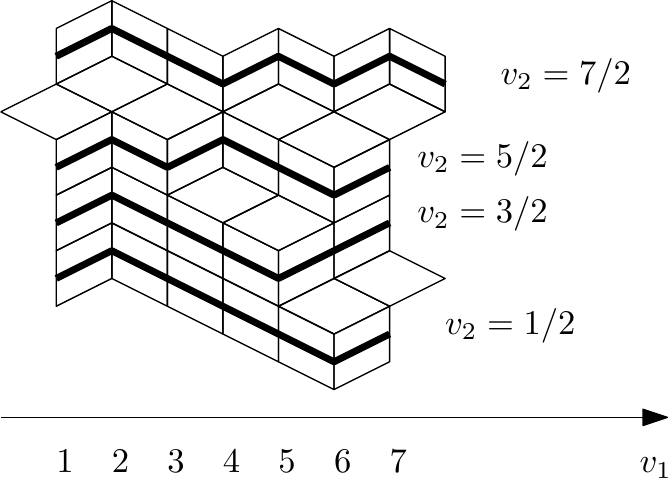}
\caption{}
\label{fig:cammini}
\end{figure}

We will call the function $v\mapsto \hat h(v)$ the ``level set
function'' of the interface $\Sigma$, in order to distinguish it from
the ``height function'' $u\mapsto h(u)$. The reason for the name is the following.
Given $v_2 \in \mathbb Z+1/2$, consider the intersection $S^{(v_2)}$
of the surface $\Sigma$ with the horizontal plane $\{(x,y,z)\in
\mathbb R^3: z=v_2\}$.  With reference to Figure \ref{fig:cammini}, each
$S^{(v_2)}$ can be viewed as a simple-random walk path
$\{S^{(v_2)}(v_1)\}_{v_1\in \mathbb Z}$ in space-time dimension $(1+1)$:
the time axis $v_1$ is horizontal, $S^{(v_2)}(v_1)\in \mathbb Z$ and
$S^{(v_2)}(v_1)-S^{(v_2)}(v_1+1)\in \{-1,+1\}$. Moreover, these lines are
mutually non-intersecting: $S^{(v_2+1)}(v_1)\ge S^{(v_2)}(v_1)$.
It is easy to check that, modulo a global additive constant independent of $v=(v_1,v_2)$, one has 
\begin{eqnarray}
  \label{eq:hatacca}
\hat h(v)= S^{(v_2)}(v_1).
\end{eqnarray}


Next, we define the analogue (in this new parametrization of the
surface) of the domain $U_L\subset \mathcal T$. Given $U_L$ as in
Section \ref{sec:monotonesurf} and $\eta\in\Omega_{U_L}$, let
$\Sigma_L$ be the monotone surface whose $\Pi_{111}$ projection is
$\eta$. Let $\eta'$ and $\Sigma'$ denote an arbitrary extension of
$\eta$ to a tiling of the whole plane $\mathcal T$ and the
corresponding monotone surface (see discussion just 
before Assumption \ref{assUL}). We let
\begin{eqnarray}
  V_L:=\Pi_{110}(\Sigma_L) 
\end{eqnarray}
and we note that $V_L$ is independent of the choice of $\eta\in
\Omega_{U_L}$. Actually, on $\mathbb
R^2\setminus V_L$ the function $v\mapsto \hat h(v)$ depends only on the arbitrary choice of $\eta'$ outside of $U_L$. 

The following is equivalent to Assumption \ref{assUL} and is actually
a rephrasing of it: 
\begin{Proposition}\label{prop:level_set_assumption}
  As $L\to\infty$, $(1/L)V_L$ tends in Hausdorff distance to a bounded
  simply connected closed domain $V\subset \mathbb R^2$ with smooth
  boundary. There is a continuous function $\hat \phi$ on $\partial V$
  and if $v^{(L)}\in \mathbb R^2\setminus V_L$ is such that
  $v^{(L)}/L\to v\in \partial V$ as $L\to\infty$ then 
  \begin{eqnarray}
\frac1L \hat h(v^{(L)})\to\hat\phi(v).
  \end{eqnarray}
Moreover, there exists a $C^1$ function $\hat \psi:V\mapsto \mathbb R$ such that 
$\hat\psi=\hat\phi$ on $\partial V$ and 
\begin{eqnarray}
\label{eq:nmass2}
 \nabla\hat\psi(v)\in (-1,1)\times (0,\infty) \quad \text{ for every } \quad v\in V.
\end{eqnarray}
\end{Proposition}
All these claims follow from Assumption \ref{assUL}  and the change of variable formulas in
Section \ref{sec:chofv} below; in particular,
Eq. \eqref{eq:sbizzarro2} shows that
\eqref{eq:nmass2} is equivalent to the non-extremality condition
\eqref{eq:nonmassimale}.

\subsubsection{Hydrodynamic limit for the level set function}

\label{sec:st2}
Let for  $v\in V$ and $t\ge0$
\begin{eqnarray}
  \hat H(v,t)=\frac1L \hat h_{\eta(t L^2)}( v L).
\end{eqnarray}
Recall that we assume that the initial condition of the dynamics
satisfies \eqref{eq:iniziale}. As in Proposition
\ref{prop:level_set_assumption}, in terms of the ``level set
function'' this implies that there exists a smooth
$\hat\psi_0:v\in V\mapsto \hat\psi_0(v)\in \mathbb R$ satisfying
\eqref{eq:nmass2} such that
\begin{eqnarray}
  \label{eq:inizialebis}
  \lim_{L\to\infty} \hat H(v,0)=\hat\psi_0(v), \quad\forall v\in V.
\end{eqnarray}
The conjectural existence of a hydrodynamic limit means existence of a function
$\psi(\cdot,\cdot):V\times \mathbb R^+\mapsto \mathbb R$, such that for
every $\epsilon>0$,
\begin{eqnarray}
  \label{eq:generale2}
  \mathbb P(\exists v\in V: |\hat H(v,t)-\hat\psi(v,t)|\ge \epsilon)\stackrel{L\to\infty}\to 0.
\end{eqnarray}
Under a (reasonable) assumption of local equilibrium, 
we find (see Section \ref{sec:localeq}) that $\hat\psi$ has to satisfy the PDE
\begin{eqnarray}
  \label{eq:PDE2}
\left\{
  \begin{array}{l}
  \partial_t \hat \psi=\frac12\partial^2_{v_1}\hat \psi-\frac14\partial_{v_2}\left[
(2+\partial_{v_2} \hat \psi)\pi_{\rho(\nabla \hat \psi)}\left[(|n(b^{(0,0)})|-1)1_{X(0,0)}1_{n(b^{(0,0)})\ne 0}
\right]
\right]    \\
\hat \psi(v,t)=\hat \psi_0(v) \quad\quad\quad\quad \text{ if } t=0 \text{ or if } v\in \partial V
  \end{array}
\right.
\end{eqnarray}
where $\rho(\nabla \hat \psi)=(\rho_1(\nabla \hat
\psi),\rho_2(\nabla\hat\psi))$ is defined as
\begin{eqnarray}
  \begin{aligned}
  \label{eq:sbizzarro}
\rho_1(s_1,s_2)=\frac{1-s_1}{2+s_2},\quad\quad
\rho_2(s_1,s_2)=\frac{1+s_1}{2+s_2}.
  \end{aligned}
\end{eqnarray}

We will prove in Section \ref{sec:K}:
\begin{Proposition}
\label{lemma:piX}
 The following identities hold:
  \begin{eqnarray}
    \label{eq:piX}
    \pi_\rho(1_{X(0,0)}1_{n(b^{(0,0)})\ne 0})=2[\rho_1 \rho_2+(1-\rho_1-\rho_2)V(\rho)]
  \end{eqnarray}
and
\begin{eqnarray}
\label{eq:orr}
  \pi_{\rho}\left[(|n(b^{(0,0)})|-1)1_{X(0,0)}1_{n(b^{(0,0)})\ne 0}
\right]=  2[-\rho_1 \rho_2+(\rho_1+\rho_2)V(\rho)].
\end{eqnarray}
\end{Proposition}
Equation \eqref{eq:PDE2} then becomes
\begin{eqnarray}
  \label{eq:limeq2}
  \partial_t \hat \psi=\frac12 \partial^2_{v_1} \hat \psi+\left.[
\partial_{ s_1} G(s)    \,\partial^2_{v_1 v_2}\hat \psi  +  \partial_{s_2} G(s)\,\partial^2_{v_2} \hat \psi   ]\right|_{s=( s_1, s_2)=\nabla \hat \psi}
\end{eqnarray}
where
\begin{eqnarray}
  \label{eq:F}
  G(s)=(1+s_2/2)\left[
\rho_1(s_1, s_2)\rho_2(s_1,s_2)-(\rho_1(s_1,s_2)+\rho_2( s_1,s_2))V(\rho(s_1, s_2))
\right]
\end{eqnarray}
and $\rho(\cdot,\cdot)$ is defined in \eqref{eq:sbizzarro}.

The derivatives of $G$ are a combination of trigonometric functions,
but 
it is best to express them in terms of the surface tension. 
Given Theorem \ref{th:CKP}, it is easy to deduce (see Section
\ref{sec:chofv}) that, given 
  $\hat\psi:V\mapsto\mathbb R$ satisfying \eqref{eq:nmass2}, one has
  \begin{eqnarray}
    \label{eq:KLP2}
  \lim_{\delta\to0}\lim_{L\to\infty}\ln   \left|\{\eta\in \Omega_{U_L}:\sup_{v\in V}|L^{-1}\hat h_\eta(vL)- \hat\psi(v)|\le \delta\}\right|=  -\int_V \hat\sigma (\nabla \hat\psi) dv
  \end{eqnarray}
where 
\begin{eqnarray}
  \label{eq:hatsigma}
\hat \sigma(s)=\hat \sigma(s_1,s_2)=\left(1+\frac{ s_2  }2\right)\sigma(\rho( s_1, s_2)).
\end{eqnarray}
Then, one can check\footnote{The right- and left-hand sides of  \eqref{eq:identita} are rather complicated trigonometric functions, and  we found no better way to check the equalities than with the help of symbolic computation using Mathematica.} from 
\eqref{eq:hatsigma} and \eqref{eq:F} that the following simple relations hold:
\begin{gather}
  \label{eq:identita}
  \partial_{s_1} G(s)=\frac{ \hat \sigma_{1,2}( s)}
{ \hat \sigma_{1,1}(s)},
\quad
\partial_{s_2} G(s)=\frac12\frac{ \hat \sigma_{2,2}(s)}
{ \hat \sigma_{1,1}(s)}
\end{gather}
with  $\hat\sigma_{i,j}:=\partial^2\hat\sigma/\partial {s_i}\partial {s_j}$.
Altogether, the conjectural hydrodynamic limit equation is 
 \begin{eqnarray}
   \label{eq:PDEter}
 \left\{
   \begin{array}{l}
   \partial_t \hat\psi=\frac12 \partial^2_{v_1} \hat \psi+\frac{\hat \sigma_{1,2}(\nabla\hat\psi)}{\hat \sigma_{1,1}(\nabla\hat\psi)}\partial^2_{v_1 v_2}\hat\psi + \frac12\frac{\hat \sigma_{2,2}(\nabla\hat\psi)}{\hat \sigma_{1,1}(\nabla\hat\psi)}\partial^2_{v_2}\hat \psi\\   
  \hat\psi(v,t)=\hat\psi_0(v) \quad\quad\quad \text{ if } t=0 \text{ or if } v\in \partial  V.
   \end{array}
 \right.
 \end{eqnarray}

\begin{Remark}
\label{rem:esotica}
  From \eqref{eq:identita} it follows that  the vector field 
  \begin{eqnarray}
    \label{eq:curlo}
s\in (-1,1)\times \mathbb R^+ \mapsto \left(\frac{\hat \sigma_{1,2}(s)}{\hat \sigma_{1,1}(s)},\frac12\frac{\hat \sigma_{2,2}(s)}{\hat \sigma_{1,1}(s)}\right)
      \end{eqnarray}
has zero curl, which of course could also be checked directly from the
definition of $\hat \sigma$. 
This is reminiscent of a  surprising identity discovered in \cite{KOS}. Namely, one easily deduces from   \eqref{eq:3} that
\begin{eqnarray}
  \label{eq:KOS}
  \det\left(
    \begin{array}{ll}
      \sigma_{1,1} & \sigma_{1,2}\\
\sigma_{1,2}& \sigma_{2,2}
    \end{array}
\right)=\pi^2.
\end{eqnarray}
(Actually, identity \eqref{eq:KOS} holds more generally for the dimer
model on any infinite, translation invariant, planar bipartite lattice
\cite[Sec. 5.3.3]{KOS} and it is deduced from algebraic properties of so-called
spectral curves \cite[Sec. 3.2.3]{KOS} of the dimer model.)
Taking the derivative of \eqref{eq:KOS} w.r.t. $\rho_1 $ one obtains
\begin{eqnarray}
  \sigma_{1,1,1}\sigma_{2,2}+\sigma_{1,1}\sigma_{1,2,2}=2 \sigma_{1,2}\sigma_{1,1,2}
\end{eqnarray}
with of course $\sigma_{a,b,c}=\partial_{\rho_a}\partial_{\rho_b}\partial_{\rho_c}\sigma(\rho)$. In turn, this is immediately seen to be  equivalent to the vanishing of the curl of the vector field
\begin{eqnarray}
  \label{eq:curlino}
\rho\in \mathbb T\mapsto \left(\frac{\sigma_{1,2}(\rho)}{\sigma_{1,1}(\rho)},\frac12\frac{ \sigma_{2,2}(\rho)}{ \sigma_{1,1}(\rho)}\right) .
\end{eqnarray}
Equation \eqref{eq:curlino} is strikingly similar to \eqref{eq:curlo}
but it is not obvious how to relate via simple algebra the vanishing of the curl of the
former with that of the latter. The probabilistic
meaning and the relation between these two identities deserves to be
further explored.

For more general random surface models, there is no reason
to expect that the l.h.s. of \eqref{eq:KOS} is independent of the
slope (convexity of the surface tension says just that the determinant
is positive) but in \cite{AAY} it was argued that \eqref{eq:KOS}
should hold at slopes corresponding to cusps of $\sigma$.  
\end{Remark}



A couple of comments are in order. First, as it should, equations
\eqref{eq:PDEter} and \eqref{eq:PDE} (together with \eqref{eq:4}) are
actually the same PDE, as can be seen with a suitable change of
variables (see Section
\ref{sec:chofv} for details).
A second remark is that in the $(v_1,v_2)$ coordinate system one
obtains a surprisingly simple relation between interface mobility and
surface tension. Namely, rewriting (in
analogy with \eqref{eq:PDE}) the PDE satisfied by $\hat \psi$ as
 \begin{eqnarray}
   \label{eq:PDEbis}
   \partial_t \hat\psi=\hat\mu(\nabla \hat\psi)\sum_{i,j=1,2}\hat\sigma_{i,j}(\nabla \hat\psi)\frac{\partial^2}{\partial_{v_i}\partial_{v_j}}\hat\psi
 \end{eqnarray}
and comparing with \eqref{eq:PDEter} one sees that the mobility coefficient is
\begin{eqnarray}
  \label{eq:muhat}
  \hat\mu(\nabla\hat \psi)=\frac1{2\hat\sigma_{1,1}(\nabla\hat \psi)}.
\end{eqnarray}
It would be interesting to understand whether there is any thermodynamic explanation for such relation.


Finally, from 
\eqref{eq:limeq2} one sees that the equation for $\hat\psi$ is also given by  
\begin{eqnarray}
\label{eq:PDEdiv2}
  \partial_t\hat \psi={\rm div}((1/2)\partial_{v_1}\hat \psi,G(\nabla \hat \psi)),
\end{eqnarray}
from which one finds again that the time derivative of the total volume is a boundary term:
\begin{eqnarray}
  \frac d{dt}\int_V \hat\psi(v,t)=\int_{\partial V}((1/2)\partial_{v_1}\hat\psi, G(\nabla \hat\psi))\cdot n\,d\gamma.
\end{eqnarray}
Take now two initial conditions with $\hat\psi^{(1)}(v)\ge
\hat\psi^{(2)}(v)$ for every $v\in V$ (or equivalently
$\psi^{(1)}(u)\ge \psi^{(2)}(u)$ for every $u\in U$) and coinciding on $\partial V$. We want to show
that
 \begin{eqnarray}
   \label{eq:drift2}
   \frac d{dt}\int_V (\hat\psi^{(1)}(v,t)-\hat\psi^{(2)}(v,t))dv\le 0.
 \end{eqnarray}
 As we already observed after \eqref{eq:driftnegativo}, the property
 $\psi^{(1)}\ge \psi^{(2)}$ is preserved at later times, by the
 comparison principle for parabolic PDEs. One deduces
\begin{eqnarray}
  \begin{aligned}
  \label{eq:dn1}
\frac d{dt}\int_V (\hat\psi^{(1)}-\hat\psi^{(2)})dv
= \int_{\partial V}((1/2)\partial_{v_1}(\hat\psi^{(1)}-\hat\psi^{(2)}),G(\nabla\hat\psi^{(1)})-G(\nabla\hat\psi^{(2)}))\cdot n\; d\gamma\\
= \frac12\int_{\partial V}n_1\partial_{v_1}(\hat\psi^{(1)}-\hat\psi^{(2)})\; d\gamma+
\int_{\partial V}n_2(G(\nabla\hat\psi^{(1)})-G(\nabla\hat\psi^{(2)}))\; d\gamma,
  \end{aligned}
\end{eqnarray}
with $n_i$ the $i^{th}$ component of the exterior normal vector $n$.
Since $\hat \psi^{(1)}\ge \hat \psi^{(2)}$ and they coincide on
$\partial V$, for $v_0\in \partial V$ we have that\footnote{Recall
  that we are assuming that $\hat\psi(\cdot,t)$ is $C^1$ in $V$ for
  all times, see discussion at the beginning of Section
  \ref{sec:vd}. }
\begin{eqnarray}
  \label{eq:segni}
  {\rm sign}\left[
\lim_{v\to v_0} \partial_{v_i}(\hat\psi^{(1)}-\hat\psi^{(2)})\right]=-{\rm sign} \left[n_i(v_0)\right],\quad i=1,2,
\end{eqnarray}
whenever $n_i(v_0)\ne 0$ (here, $n(v_0)$ is the exterior normal at $v_0$).  To show \eqref{eq:drift2} (and actually that the integrand of \eqref{eq:dn1} itself is point-wise negative) it
suffices to show that
\begin{eqnarray}
  \label{eq:segni2}
  {\rm sign}\left[ \lim_{v\to v_0} \partial_{v_2}(\hat\psi^{(1)}-\hat\psi^{(2)})\right]= {\rm sign}\left[
\lim_{v\to v_0}(G(\nabla\hat\psi^{(1)})-G(\nabla\hat\psi^{(2)}))
\right]
\end{eqnarray}
whenever $n_2(v_0)\ne 0$.
Take the limit ${v\to v_0}$ in direction $v_2$. Then, 
the claim follows if we know that $G(s_1,s_2)$ is increasing in its second argument. In turn, this follows from the second of \eqref{eq:identita} plus convexity of the surface tension $\hat\sigma$.

From formulas \eqref{eq:11} and \eqref{eq:15} of Section \ref{sec:chofv} one sees that\footnote{The factor $1/2$ in Eq. \eqref{eq:drift1} can be intuitively understood from the microscopic dynamics, recalling Remark \ref{rem:duevolumi}.}
\begin{eqnarray}
  \label{eq:drift1}
 \frac d{dt}\int_U \psi(u,t) du=\frac12 \frac d{dt}\int_V \hat\psi(v,t) dv.
\end{eqnarray}
Then, \eqref{eq:driftnegativo} follows.

\section{Mobility and linear response theory}
\label{sec:GK}

Let us briefly recall the method of \cite{Spohn} for the computation
of the mobility. The idea is to consider the dynamics in a finite
volume $\Lambda$ of diameter $O(L)$ and to 
modify the rates through a ``magnetic field'' $B>0$ such that the
invariant measure $\pi^{(B)}_\Lambda$ is tilted by the exponential of $B$ times the volume
(the volume being the sum of the heights in $\Lambda$) w.r.t. the
$B=0$ situation. In presence of magnetic field,
the free energy \eqref{eq:6} becomes 
$F(\psi)= \int_U \sigma(\nabla \psi)du-B\int_U \psi du$, where the surface tension $\sigma$ is as in Theorem \ref{th:CKP}. The transition rates are changed in an easy way by $B$:
If for $B=0$   the rate at which the vertical position of particle $b$ changes from $n(b)$ (the present position in $\eta$) to $n(b)+y$, $y\in \mathbb Z,y\ne0$ is denoted $c_{b,y}(\eta)$, 
more explicitly (recalling Definitions \ref{def1} and \ref{def2}) 
\begin{eqnarray}
c^I_{b,y}(\eta)=\frac1{2|y|}1_{\{ n(b)+y\in I(b)\}}  
\end{eqnarray}
and
\begin{eqnarray}
  \label{eq:cLE}
  c^{II}_{b,y}(\eta)=\frac1{|I(b)|}1_{\{ n(b)+y\in I(b)\}},
\end{eqnarray}
for dynamics I and II respectively, then for $B\ne0$ one has
\begin{eqnarray}
  \label{eq:cB}
c^{B}_{b,y}(\eta)=c_{b,y}(\eta)\exp(By/2).  
\end{eqnarray}

The heuristics that leads to \eqref{eq:7} should still apply so the hydrodynamic limit should follow the equation
\begin{multline}
\partial_t \psi = - \mu(\nabla\psi) \frac{\delta}{\delta \psi} \left( \int_U \sigma(\nabla \psi)du - B \int_U \psi du \right) \\= \mu(\nabla \psi)\sum_{i,j=1,2}\sigma_{i,j}(\nabla \psi)\frac{\partial^2}{\partial_{u_i}\partial_{u_j}}\psi + \mu(\nabla \psi) B.
\end{multline}
If the configuration $\psi$ is flat (i.e. has  constant gradient
$\rho$ on $U$), one
has then
\begin{eqnarray}
  \label{eq:8}
  \partial_t \psi=B\mu(\rho)
\end{eqnarray}
so the $B=0$ value of the mobility is obtained (within linear response theory) as the
derivative of the interface velocity w.r.t. $B$, at $B=0$. The system is
started from the $B=0$ reversible measure $\pi_\Lambda$ on $\Lambda$,
with boundary conditions chosen such that the law of the
gradients is translation-invariant and corresponds to the desired slope
$\rho$. Then, the box side $L$ is taken
to $+\infty$ at fixed time (so that boundary effects do not prevent the law of
the gradients from remaining translation invariant at positive time)
and finally time is taken to $+\infty$, in order to ensure that the
system has reached its asymptotic velocity. Note that if time were
taken to infinity first, then the asymptotic velocity  would be $0$
since the system would reach its new invariant measure $\pi^{(B)}_\Lambda$.

Let us implement this scheme in our case. For technical convenience,
we take $\Lambda\subset \mathcal T$ to be $L$-periodic both in
direction $e_1$ and in direction $e_2$. Note that
any tiling of 
$\Lambda$ contains $L^2$ lozenges. In order to impose the
tilt $\rho$, we let $\pi_\Lambda:=\pi_{\Lambda, \rho^{(L)}}$ denote the uniform measure on lozenge
tilings of $\Lambda$ such that the fraction of north-west
(resp. north-east) oriented
lozenges  is  $\rho^{(L)}_1$ (resp. $\rho^{(L)}_2$), with
$\rho_i^{(L)}\to \rho_i$ as $L\to\infty$. It is then known \cite{KOS} that
$\pi_\Lambda$ converges to $\pi_{\rho}$, in the sense that
averages of local functions converge.

The drift of height per unit site (i.e. the interface velocity) at time $t$ is
\begin{gather}
\label{eq:precise}
\frac1{L^2}  \sum_{b,y}\mathbb E_{\pi_\Lambda}\left[ c^B_{b,y}(\eta(t))y\right]=\frac1{L^2}  \sum_{\eta,\eta'} \sum_{b,y}\pi_\Lambda(\eta)e^{t \mathcal L^B}(\eta,\eta')c^B_{b,y}(\eta')y
\end{gather}
with the average on the law of the process at time $t$ started from
the $\pi_\Lambda$ and the sum over $b$ running over particles in
$\Lambda$ (there are exactly $L^2(1-\rho_1^{(L)}-\rho_2^{(L)})$ of them).
Here $\mathcal L^B$ is the generator of the process with rates $c^B_{b,y}$.

Recall that we want to take the $B$ derivative of \eqref{eq:precise} at $B=0$, then the
limit $L\to \infty$ and finally the limit $t\to\infty$ to obtain the
mobility $\mu(\rho)$.
Taking the derivative and using 
\begin{eqnarray}
  \label{eq:masu}
\left.\frac d{dB}c^B_{b,y}(\eta)\right|_{B=0}=c_{b,y}(\eta)\frac y2  
\end{eqnarray}
together with invariance of $\pi_\Lambda$ for the  $B=0$ process,
one gets
\begin{gather}
\label{deriv}
\frac12\frac1{L^2}  \sum_{b,y}  \sum_\eta \pi_\Lambda(\eta)c_{b,y}(\eta)y^2
+\frac1{L^2} \sum_{\eta,b,y,\eta',\eta''}\int_0^t 
\pi_\Lambda(\eta) \mathcal L'(\eta,\eta')e^{\tau\mathcal L}(\eta',\eta'')c_{b,y}(\eta'')y \, d\tau
\end{gather}
where $\mathcal{L}'(\eta,\eta')=\partial_B \mathcal L(\eta,\eta')|_{B=0}$. 
We have
\begin{gather}
  \sum_{\eta,\eta',\eta''}\pi_\Lambda(\eta) \mathcal L'(\eta,\eta')e^{\tau\mathcal L}(\eta',\eta'')c_{b,y}(\eta'')y\\=
\label{somma}
  \sum_{\eta,\eta'}\pi_\Lambda(\eta) \mathcal L'(\eta,\eta)e^{\tau\mathcal L}(\eta,\eta')c_{b,y}(\eta')y
+ \sum_{\eta,\eta',\eta'':\eta'\ne\eta}\pi_\Lambda(\eta) \mathcal L'(\eta,\eta')e^{\tau\mathcal L}(\eta',\eta'')c_{b,y}(\eta'')y.
\end{gather}
Note that $\mathcal L'(\eta, \eta')$ (with $\eta\ne\eta'$) is non-zero
only if one can go from $\eta$ to $\eta'$ in a single move (say $b', z$)
and in that case $\mathcal L'(\eta, \eta') = c_{b', z}(\eta)
\frac{z}{2}$. Then, the last
sum becomes
\begin{gather}
  \sum_{\eta,\eta'',b',z}\pi_{\Lambda}(\eta)c_{b',z}(\eta)\frac z2 e^{\tau \mathcal L}(\eta^{b',z},\eta'')c_{b,y}(\eta'')y
\end{gather}
with $\eta^{b',z}$ the configuration where $b'$ has been moved by $z$.
Changing names of variables this is also equal to
\begin{gather}
  \sum_{b',z}\sum_{\eta,\eta'}\pi_{\Lambda}(\eta^{b',-z})c_{b',z}(\eta^{b',-z})\frac z2 e^{\tau \mathcal L}(\eta,\eta')c_{b,y}(\eta')y
\\=-\frac12   \sum_{b',z}\sum_{\eta,\eta'}\pi_{\Lambda}(\eta^{b',z})c_{b',-z}(\eta^{b',z}) z e^{\tau \mathcal L}(\eta,\eta')c_{b,y}(\eta')y.
\end{gather}
Using detailed balance, $\pi_{\Lambda}(\eta^{b',z})c_{b',-z}(\eta^{b',z})=\pi_\Lambda(\eta)c_{b',z}(\eta)$.
Therefore \eqref{somma} is equal to 
\begin{gather}
\label{somma2}
  \sum_{\eta,\eta'}\pi_\Lambda(\eta) [\mathcal L'(\eta,\eta)-\sum_{\eta'':\eta''\ne\eta}\mathcal L'(\eta,\eta'')]e^{\tau\mathcal L}(\eta,\eta')c_{b,y}(\eta')y  .
\end{gather}
Since $\sum_{\eta''} \mathcal L(\eta,\eta'')=0$, the same is true for $\mathcal L'$ and \eqref{somma2} gives
\begin{gather}
 -2   \sum_{\eta,\eta',\eta'':\eta''\ne\eta}\pi_\Lambda(\eta) \mathcal L'(\eta,\eta'')e^{\tau\mathcal L}(\eta,\eta')c_{b,y}(\eta')y.
\end{gather}
Altogether, using also \eqref{eq:masu}, \eqref{deriv}
gives
\begin{gather}
\label{mob1}
 \frac1{2L^2}\pi_\Lambda[\sum_{b,y}c_{b,y}(\eta)y^2]-\frac1{L^2}\sum_{b,y,b',z}\int_0^t  \mathbb P_{\pi_\Lambda}\left[
c_{b',z}(\eta(0))z \;c_{b,y}(\eta(\tau))y\right] d\tau.
\end{gather}
Finally, letting first $L\to\infty$ and then $t\to\infty$, the mobility provided by  the linear response theory is 
\begin{gather}
  \label{eq:muu}
  \mu(\rho)=\lim_{L\to\infty}
  \frac1{2L^2}\pi_\Lambda[\sum_{b,y}c_{b,y}(\eta)y^2]
-\int_0^\infty\lim_{L\to\infty}\frac1{L^2}\sum_{b,y,b',z} \mathbb P_{\pi_\Lambda}\left[
c_{b',z}(\eta(0))z \;c_{b,y}(\eta(\tau))y\right]d\tau.
\end{gather}
A simple computation shows that, both for  dynamics I and II,
\begin{eqnarray}
  \label{eq:summa}
\pi_\Lambda[\sum_y  c_{b,y}(\eta)y^2|n^+(b),n^-(b)]=\frac1{6}(|I(b)|^2-1)=\frac1{6}(|I(b)|+1)(|I(b)|-1),
\end{eqnarray}
where we recall that $n^+(b),n^-(b)$ are the highest/lowest vertical position particle $b$ can take, given the positions of the others, and $I(b)=\{n^-(b),\dots,n^+(b)\}$. 
Let us see why. The point is that, under the uniform measure $\pi_\Lambda$ and conditionally on $n^\pm(b)$, the position 
$n(b)$ of particle $b$ is uniform in $I(b)$. Then, for dynamics I we have
\begin{gather}
  \pi_\Lambda[\sum_y  c_{b,y}(\eta)y^2|n^+(b),n^-(b)]=\frac12  \pi_\Lambda\left[\left.\sum_{y: n^-(b)\le n(b)+y\le n^+(b)} \!\!\!\!\!\!\!\!\!\!\!\!\!\!\!|y| \;\;\;\;\;\;\right|n^+(b),n^-(b)\right]\\=
\frac12 |I(b)| \mathbb E|U-V|,
\end{gather}
with $U,V$ two independent random variables, uniformly distributed in $I(b)$. An immediate calculation shows that
\[\mathbb E|U-V|=\frac13\frac{|I(b)|^2-1}{|I(b)|}
\]
and 
 \eqref{eq:summa} follows. For dynamics II one has instead
 \begin{gather}
     \pi_\Lambda[\sum_y  c_{b,y}(\eta)y^2|n^+(b),n^-(b)]=\frac1{|I(b)|}   \pi_\Lambda\left[\left.\sum_{y:n^-(b)\le n(b)+y\le n^+(b)} \!\!\!\!\!\!\!\!\!\!\!\!\!\!\!y^2 \;\;\;\;\;\;\right|n^+(b),n^-(b)\right]\\
=\mathbb E|U-V|^2,
 \end{gather}
 with $U,V$ uniform on $I(b)$ as before. A simple computation
 gives \[\mathbb E|U-V|^2=\frac16(|I(b)|^2-1)\]
 and again \eqref{eq:summa} follows.

Let us resume from Eqs. \eqref{mob1} and \eqref{eq:muu} the computation of the mobility.  For every particle $b$, there exist exactly $|I(b)|$ sites $u\in
\Lambda$ such that $b^{(u)}=b$ (recall that $b^{(u)}$ was defined just
before \eqref{eq:mu}) and exactly one of those sites satisfies
$(u_1+u_2)/2=n(b)=n(b^{(u)})$.  Therefore, the first term in \eqref{mob1} is
\begin{eqnarray}
  \label{eq:mob2}
  \frac1{12L^2}\pi_\Lambda[\sum_{u\in\Lambda}(|I(b^{(u)})|+1)1_{X(u)} 1_{n(b^{(u)})\ne (u_1+u_2)/2}],
\end{eqnarray}
where the indicator function that $n(b^{(u)})\ne (u_1+u_2)/2$ reconstructs the factor $(|I(b)|-1)$ in \eqref{eq:summa} and the indicator $1_{X(u)}$ excludes the sites $u$ for which $b^{(u)}$ is not defined.
By translation invariance and convergence of $\pi_\Lambda$ to
$\pi_{\rho}$ as  $L\to\infty$, this converges to\footnote{strictly speaking, $I(b^{(0,0)})$ and $1_{X((0,0))}$ are
    not local observables, so proving convergence requires some (doable)
    technical work. Since anyway linear response theory involves other
non-rigorous steps that are much harder to justify,
we do not give details on this point. }
\begin{eqnarray}
  \label{eq:mob3}
\frac1{12}\pi_{\rho}( (|I(b^{(0,0)})|+1) 1_{X((0,0))}1_{n(b^{(0,0)})\ne 0}).
\end{eqnarray}
Now let us look at the second term in \eqref{eq:muu}. Both for   dynamics I and II one has, by direct computation,
\begin{eqnarray}
  \sum_y c_{b,y}(\eta)y=\frac12[n^+(b)-2n(b)+n^-(b)]=\frac12[|I^+(b)|-|I^-(b)|]
\end{eqnarray}
where $I^+(b)=\{k\in I(b):k>n(b)\}$ and $I^-(b)=\{k\in I(b):k<n(b)\}$ (i.e. the set of positions available for $b$ above/below the present position). It was proven in \cite[Sec. 4]{Toninelli2+1} that
\begin{eqnarray}
  \sum_b [|I^+(b)|-|I^-(b)|]=0
\end{eqnarray}
for any configuration. (Here it is important that we are working on
the torus $\Lambda$ to avoid boundary terms).
Altogether, \eqref{eq:mu} follows.

\section{Derivation of the hydrodynamic equation under assumption of local equilibrium}
\label{sec:LE}

\subsection{Changes of variables}

\label{sec:chofv}

Given a function $\psi:u\in U\mapsto \psi(u)$ satisfying \eqref{eq:nonmassimale}, call
$\Gamma\subset \mathbb R^3$ the surface whose height function is
$\psi$: the point $p=(x,y,z)\in \Gamma$ whose $\Pi_{111}$ projection is  $u\in U$ (i.e. such that $p=(a,a,a)+(u_1,u_2,0)$ for some $a$) has vertical coordinate $z=a=-\psi(u)$.
Also  let $\hat\psi: v\in V\mapsto \hat\psi(v)$ be the
corresponding level set function: the point $p=(x,y,z)\in \Gamma$ whose $\Pi_{110}$ projection is $v$ (i.e. 
$p=(b,b,0)+(-v_1/2,v_1/2,v_2)$ for some $b$)
satisfies 
$- p\cdot(e_x+e_y) = -2b=\hat \psi(v)$, see \eqref{eq:hatac}.
Conversely, if $(x,y,z)\in \Gamma$ then its $(u_1,u_2)$ and $(v_1,v_2)$ coordinates satisfy
\begin{eqnarray}
  \label{eq:sgc1}
  u_1=x-z,u_2=y-z, \psi(u)=-z
\end{eqnarray}
and 
\begin{eqnarray}
  \label{eq:sgc2}
  v_1=y-x, v_2=z, \hat \psi(v)=-x-y.
\end{eqnarray}

Given $u\in U$, let $v(u)\in V$ be the
 $\Pi_{110}$ projection of the unique point $p\in\Gamma$ whose
$\Pi_{111}$ projection is $u$. The function $u:v\in V\mapsto u(v)$ will
denote the inverse\footnote{The function $v(\cdot)$ is
  bijective thanks to \eqref{eq:nonmassimale}: bijectivity could be lost if one allowed $\nabla\psi$ to take the value $(0,0)$, since horizontal regions of $\Gamma$ with normal vector $(0,0,1)$ project into segments of the plane $P_{110}$. } of $v(\cdot)$. 
By simple geometric considerations (see the end of this subsection),  one finds that
 the area element $du=du_1du_2$ equals  
 \begin{eqnarray}
   \label{eq:15}
   du=\frac{dv}{\partial_{u_1}\psi(u(v))+\partial_{u_2}\psi(u(v))}.
 \end{eqnarray}
 One can check this also with a microscopic argument: a finite portion
 $\Sigma'\subset \Sigma$ composed of $N$ square faces, once projected
 on $P_{111}$, contains $N$ lozenges and has area $N$ in the $u_1,u_2$
 coordinates. Assume that $N\rho_i $ lozenges are of type $i$.
When one projects $\Sigma'$ on $P_{110}$ instead, horizontal faces are
mapped to zero-area segments and the other faces to area-1 squares. 
Altogether, the area of the $P_{110}$ projection is $N(\rho_1+\rho_2)$
and \eqref{eq:15} follows from the identification between
$\partial_{u_i}\psi, i=1,2$ and $\rho_i$.

Moreover, one has
\begin{eqnarray}
  \label{eq:sbizzarro2}
  \begin{aligned}
\partial_{u_1}\psi(u)=\frac{1-\partial_{v_1}\hat\psi(v(u))}{2+\partial_{v_2}\hat \psi(v(u))}=\rho_1(\nabla\hat \psi(v(u)))\\
\partial_{u_2}\psi(u)=\frac{1+\partial_{v_1}\hat \psi(v(u))}{2+\partial_{v_2}\hat \psi(v(u))}=\rho_2(\nabla \hat \psi(v(u))),    
  \end{aligned}
\end{eqnarray}
where $\rho(\cdot)$ was already defined in \eqref{eq:sbizzarro}.
To see this, call $n=(n_1,n_2,n_3)$ the normal vector of the surface $\Gamma$ at some point $p_0=(x_0,y_0,z_0)$, normalized so that $n_1+n_2+n_3=1$. Let us assume w.l.o.g. that $p_0=(0,0,0)$.
At first order around $p_0$, the coordinates of a point $p=(x,y,z)\in\Gamma$ verify
\begin{eqnarray}
  \label{eq:sgc3}
x n_1+y n_2+z n_3=O(|p|^2).  
\end{eqnarray}
One then immediately verifies that the height function $\psi(\cdot)$ and the level set function $\hat \psi(\cdot)$ verify
\begin{multline}
  \label{eq:upo}
  \psi(u)=n_1 u_1+n_2 u_2+O(|u|^2)\\
\hat \psi(v)=\frac{(n_2-n_1)v_1+2 n_3 v_2}{n_1+n_2}+O(|v|^2),
\end{multline}
so that 
\begin{multline}
  \label{eq:corrisp2}
  \partial_{v_1}\hat \psi(0)= \frac{n_2-n_1}{n_1+n_2}= \frac{\partial_{u_2}\psi(0)-
\partial_{u_1}\psi(0)}{\partial_{u_1}\psi(0)+\partial_{u_2}\psi(0)}\\
  \partial_{v_2}\hat \psi(0)= \frac{2n_3}{n_1+n_2}=\frac{2(1-\partial_{u_1}\psi(0)+\partial_{u_2}\psi(0))}{\partial_{u_1}\psi(0)+\partial_{u_2}\psi(0)},
\end{multline}
whose inverse is \eqref{eq:sbizzarro2}.
Note that $\rho(\cdot)$ is a bijection from $(-1,1)\times (0,\infty)$ to
the triangle $\mathbb T$ defined just after \eqref{eq:nonmassimale}.

Using \eqref{eq:sgc1} and \eqref{eq:sgc2} we see that the
infinitesimal vector $(dx, dy,dz)$ (with $dz=-(n_1 dx+n_2 dy)/n_3$ in
view of \eqref{eq:sgc3}) corresponds in the $P_{111}$ coordinates to
\[
(du_1,d u_2)=\left(\frac{n_3+n_2}{n_3}dx+  \frac{n_2}{n_3}d y,  \frac{n_1}{n_3}dx+\frac{n_2+n_3}{n_3}dy\right)
\]
and in the $P_{110}$ coordinates to
\[
(dv_1,dv_2)=\left(dy-dx,-\frac{n_1 dx+n_2 dy}{n_3}\right).
\]
Therefore, the ratio of area elements $du/dv$ equals (using $n_1+n_2+n_3=1$ and $n_1=\partial_{u_1}\psi,n_2=\partial_{u_2}\psi$ from \eqref{eq:sgc3})
\[
\frac{du}{dv}=\frac 1{n_3}\times \frac{n_3}{n_1+n_2}=\frac1{n_1+n_2}=\frac1{\partial_{u_1}\psi+\partial_{u_2}\psi}
\]
which is what we claim in \eqref{eq:15}.

\subsubsection{Relation between $\sigma$ and $\hat \sigma$}

We can now write
\begin{eqnarray}
  \label{eq:16}
  \begin{aligned}
 F(\psi)&=\int_U \sigma(\nabla \psi)du=\int_V
  \sigma(\rho(\nabla\hat\psi(v)))
  \frac{dv}{\partial_{u_1}\psi(u(v))+\partial_{u_2}\psi(u(v))}\\
&=\int_V \sigma(\rho(\nabla\hat\psi(v)))(1+\frac12\partial_{v_2}\hat
\psi(v))dv    
  \end{aligned}
\end{eqnarray}
where in the last equality we used \eqref{eq:sbizzarro2}.  This shows
that the surface tension functional can also be expressed in the
parametrization of the surface by $\hat \psi$ and that we only need to
replace $\sigma(\nabla \psi)$ by $\hat \sigma (\nabla\hat\psi)=
\sigma(\rho(\nabla\hat\psi))(1+\frac12\partial_{v_2}\hat \psi)$ as per
\eqref{eq:hatsigma}.

\subsubsection{Identity between \eqref{eq:PDE} and \eqref{eq:PDEter}}

Here we verify that the equations \eqref{eq:PDE} and \eqref{eq:PDEter} (together with
\eqref{eq:4}) for the hydrodynamic limit in the $\psi$ and $\hat \psi$ parametrization indeed describe the
same interface evolution.
Indeed, \eqref{eq:PDE} and \eqref{eq:PDEter} can be rewritten as
\begin{eqnarray}
  \label{eq:9}
  \partial_t\psi(u,t)=-\mu(\nabla \psi(u,t))\left.\frac{\delta F}{\delta \psi(u)}\right|_{\psi(\cdot):=\psi(\cdot,t)}
\end{eqnarray}
and 
\begin{eqnarray}
  \label{eq:10}
  \partial_t\hat\psi(v,t)=-\hat\mu(\nabla\hat\psi(v,t))
  \left.\frac{\delta \hat F}{\delta \hat\psi(v)}\right|_{\hat\psi(\cdot):=\hat\psi(\cdot,t)}
\end{eqnarray}
where $F(\psi)=\int_U \sigma(\nabla \psi)du$,
$\hat F(\hat\psi)=\int_V\hat\sigma(\nabla\hat\psi)dv=F(\psi)$ and
$\hat \mu$ is given in \eqref{eq:muhat}.  

One verifies by direct inspection that
\begin{eqnarray}
  \label{eq:5}
 \frac{ \mu(\nabla \psi(u,t))}{\hat \mu (s(\nabla\psi(v(u),t)))}=\frac{\partial_{u_1}\psi(u,t)+\partial_{u_2}\psi(u,t)}4,
\end{eqnarray}
where $s(\rho)$ is the inverse of the function
$\rho(s)=(\rho_1(s),\rho_2(s))$ defined in \eqref{eq:sbizzarro} or \eqref{eq:sbizzarro2}.
Next, we will prove that 
\begin{eqnarray}
  \label{eq:11}
  \frac{\partial_t\psi(u,t)}{\partial_t\hat\psi(v,t)|_{v=v(u)}}=\frac{\partial_{u_1}\psi(u,t)+\partial_{u_2}\psi(u,t)}2.
\end{eqnarray}
and that 
\begin{eqnarray}
  \label{eq:12}
  \frac{\delta F}{\delta
  \psi(u)}=2
\left.\frac{\delta \hat F}{\delta\hat \psi(v)}\right|_{v=v(u)}.
\end{eqnarray}
Together with \eqref{eq:5}, these imply that \eqref{eq:9} and \eqref{eq:10} are indeed equivalent.

Let us verify \eqref{eq:11}. If $\psi(\cdot,t),\hat\psi(\cdot,t)$ are
the height function and level set function of the same evolving
surface $\Gamma_t\subset \mathbb R^3$, we want to find the relation
between $\partial_t\psi(u,t)$ and $\partial_t\hat\psi(v,t)$.
Let $p\in\Gamma_{t}$ and
let the interface normal velocity  $p$ be $V\times(n_1,n_2,n_3)$,
with $(n_1,n_2,n_3)$ the normal vector normalized as
$n_1+n_2+n_3=1$. Let the  $\Pi_{111}$ and $\Pi_{110}$ projections of $p$ have coordinates $u_0$ and $v_0$ respectively. In analogy with \eqref{eq:upo}, one has
\begin{multline}
  \label{eq:vpo}
\psi(u_0+u,t+\epsilon)=-V'\epsilon+n_1u_1+n_2u_2+O(\epsilon^2)+O(|u|^2)+C(u_0)\\
\hat\psi(v_0+v,t+\epsilon)=\frac{-2V'\epsilon+v_1(n_2-n_1)+2v_2n_3}{n_1+n_2}+O(\epsilon^2)+O(|v|^2)+C'(v_0),
\end{multline}
with $V'=(n_1^2+n_2^2+n_3^2)V$ and $C(u_0),C'(v_0)$ independent of time and of $u,v$, whence 
\begin{eqnarray}
  \label{eq:aaaga}
\frac{\partial_t \psi(u_0,t)}{\partial_t\hat \psi(v_0,t)}=\frac{n_1+n_2}2=
\frac{\partial_{u_1}\psi(u_0,t)+\partial_{u_2}\psi(u_0,t)}2  
\end{eqnarray}
and \eqref{eq:11} follows.

It remains only to check \eqref{eq:12}.
Under an infinitesimal variation $\psi\to\psi+\varepsilon \delta\psi$,
the first variation of $F$ is 
\begin{eqnarray}
  \label{eq:14}
  \int_U \delta\psi(u) \frac{\delta F}{\delta
  \psi(u)} du.
\end{eqnarray}
The same  considerations (steps \eqref{eq:vpo} and \eqref{eq:aaaga}) that lead to \eqref{eq:11} show that
\[
\delta\psi(u)=\delta\hat\psi(v(u))\frac{\partial_{u_1}\psi(u)+\partial_{u_2}\psi(u)}2.
\]
Then, together with \eqref{eq:15}, the quantity in \eqref{eq:14} can
be rewritten as 
\begin{eqnarray}
  \label{eq:17}
\frac12  \int_V \left.\frac{\delta F}{\delta\psi(u)}\right|_{u=u(v)}\delta
  \hat\psi(v) dv.
\end{eqnarray}
Since by construction $ F(\psi)=\hat F(\hat\psi)$, \eqref{eq:17}
implies \eqref{eq:12}.

\subsection{Derivation of \eqref{eq:PDE2}}
\label{sec:localeq}

We will make an assumption of local equilibrium, namely, we assume
that the law of ${\eta(L^2t)}$ in the neighborhood of any site $u_L\in
U_L$, such that $\bar u=\lim_{L\to\infty} u_L/L$ is in the interior of
$U$, approaches (when $L\to\infty$) the equilibrium measure $\pi_\rho$
with slope
$\rho=(\partial_{u_1}\psi(u,t),\partial_{u_2}\psi(u,t))=\rho(\nabla\hat\psi(v(u),t))$,
where the function $\rho(\cdot)$ is defined in
\eqref{eq:sbizzarro}. 

Let us proceed to the derivation of \eqref{eq:PDE2} under such assumption.
Let $f:V\mapsto \mathbb R$ be a $C^\infty$ test function whose compact support is at positive distance from $\partial V$.
We define for any configuration $\eta$
\begin{eqnarray}
  \label{eq:test1}
  I_f^{(L)}(\eta):=\frac1{|V_L|}\sum_{v\in V_L}f(v/L)\frac{\hat h_{\eta}(v)}L,
\end{eqnarray}
with $|V_L|=O(L^2)$ the cardinality of $V_L$. 
We will argue that 
\begin{multline}
  \label{eq:hypo2}
  \partial_t \mathbb E\left( I_f^{(L)}(\eta({L^2 t}))\right)\stackrel{L\to\infty}\to
\int_V \left[\partial^2_{v_1}f(v)\frac{\hat\psi(v,t)}2 
+\partial_{v_2}f(v)\left(1+\frac{\partial_{v_2}\hat\psi}2\right)
\pi_{\rho(\nabla\hat\psi(v,t))}(A)\right]
   dv
\end{multline}
with 
\[
A=A(\eta)= 1_{X(0,0)}1_{n(b^{(0,0)})\ne0}\frac{|n(b^{(0,0)})|-1}2\] so
that, at least in a weak sense, \eqref{eq:PDE2} must hold.  
Note also  that
in  \eqref{eq:hypo2} we are just looking at the average interface evolution. For an
actual proof of the hydrodynamic limit, one should also prove that
$I_f^{(L)}(\eta(L^2 t))-I_f^{(L)}(\eta(0))$ concentrates for $L\to\infty$ at the
time-integral of the r.h.s. of \eqref{eq:hypo2}.

The l.h.s. of \eqref{eq:hypo2} is
\begin{eqnarray}
L^2  \mathbb E\left([\mathcal L I^{(L)}_f](\eta({L^2 t}))\right),
\end{eqnarray}
with $\mathcal L$ the generator of the dynamics. We will work with 
Dynamics I (Definition \ref{def1}) and leave to the reader to check that the
same result is obtained with Dynamics II.

\begin{figure}[h]
  \includegraphics[width=4cm]{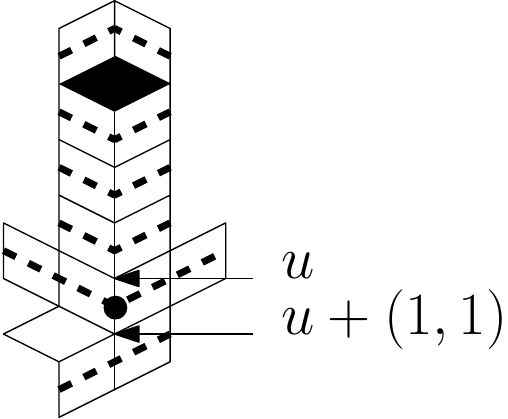}
  \caption{The edge $u,u+(1,1)$ is common to a type-1 and a type-2
    lozenge iff a level line $S^{(v_2)}$ (dashed) runs through the
    point $p $ (dot in the picture) at distance $1/2$ below $u$, and
    $S^{(v_2)}(\cdot)$ has a maximum/minimum at $v_1$ (a minimum in
    this example). Here, $(v_1,v_2)$ is the $\Pi_{110}$ projection of
    the point $p$.  The distance $|n^{(u)}-(u_1+u_2)/2|$ ($=4$ in the
    drawing) equals also $k(\eta,(v_1,v_2))$. One sees that
    $k(\eta,(v_1,v_2))$ is also the minimal value of $n\ge1$ such that
    $\hat h(v_1,v_2+n)\ne \hat h(v_1,v_2)$, or the
    number of minima of the level set function that need to be flipped in order to flip the
    minimum at $v$.}
\label{fig:ennesima}
\end{figure}

We need a couple of definitions.
Given $v\in V_L$ and a configuration $\eta$, let
\begin{eqnarray}
  \label{eq:indv}
  r(\eta,v)=1_{(\Delta_{v_1}\hat h_\eta)(v)\ne 0}=\frac12|(\Delta_{v_1}\hat h_\eta)(v)|, 
\end{eqnarray}
with $\Delta_{v_1}$ the discrete Laplacian in the direction $v_1$: $(\Delta_{v_1}\hat h_\eta)(v)=\hat h_\eta(v+(1,0))+\hat h_\eta(v-(1,0))-2
\hat h_\eta(v)$, as well as
\begin{eqnarray}
   \label{eq:epv}
\epsilon(\eta,v)   =r(\eta,v) {\rm sign}[(\Delta_{v_1}\hat h_\eta)(v)]= \frac12(\Delta_{v_1}\hat h_\eta)(v)
 \end{eqnarray}
and
\begin{eqnarray}
  \label{eq:kv}
  k(\eta,v)=r(\eta,v) \min\{n\ge1:\hat h_\eta(v+\epsilon(\eta,v)(0,n))\ne\hat
  h_\eta(v)\}.
\end{eqnarray}
 In words, $r(\eta,(v_1,v_2))$ is the indicator function that the function
 $v_1\mapsto S^{(v_2)}(v_1)$ has a local extremum at $v_1$,
 $\epsilon(\eta,v)$ is $-1 $ or $+1$ if the extremum is a maximum or
 minimum respectively and $k(\eta,v)$ is the number of minima or
 maxima that need to be flipped in order to reverse the
 minimum/maximum of  $S^{(v_2)}(\cdot)$ into a maximum/minimum, see also Figure \ref{fig:ennesima}.
 

 Recall Definition \ref{def1} of the dynamics which, with the present
 notations, reads as follows: for any $v\in V_L$ such
 that $\epsilon(\eta,v)\ne0$, with rate $1/(2k(\eta,v))$ all the
 heights $\hat h(v+(0,n)\epsilon(\eta,v)), 0\le n<k(\eta,v)$ change by
 $+2\epsilon(\eta,v)$. Therefore we have
\begin{gather}
\label{generatore}
L^2 [\mathcal L^I I^{(L)}_f](\eta)= \frac L{|V_L|}\sum_{v\in V_L}
r(\eta,v)\frac{2\epsilon(\eta,v)}{2 k(\eta,v)}
\sum_{n=0}^{k(\eta,v)-1}f\left(\frac{v_1}L,
\frac{v_2+\epsilon(\eta,v)n}L\right)\\\label{eq:somma1}
=\frac{L}{|V_L|}\sum_{v\in V_L}\frac12 [\Delta_{v_1}\hat h_\eta](v)f(v/L)\\\label{eq:somma2}+
\frac L{|V_L|}\sum_{v\in
  V_L}r(\eta,v)\frac{\epsilon(\eta,v)}{k(\eta,v)}\sum_{n=0}^{k(\eta,v)-1}\left(f\left(\frac{v_1}L,\frac{
v_2+\epsilon(\eta,v)n}L\right)-f\left(\frac vL\right)\right).
\end{gather}
The  sum in \eqref{eq:somma1} can be rewritten, via a summation by parts and using $\Delta_{v_1}f(v/L)=L^{-2}\partial^2_{v_1}f(v/L)+O(L^{-3})$, 
as 
\begin{eqnarray}
 \frac12 \frac1{|V_L|}\sum_{v\in V_L}\frac{\hat h_\eta(v)}L \partial^2_{v_1}f(v/L)+O(1/L)
\end{eqnarray}
and the error term is uniform w.r.t. $\eta$. 
There are no boundary terms in the summation by parts because $f$ is
compactly 
supported in the interior of $V$.

As for  \eqref{eq:somma2}, let us Taylor expand
\begin{eqnarray}
  \label{eq:Taylor}
 f\left(\frac{v_1}L,\frac{
v_2+\epsilon(\eta,v)n}L\right)-f\left(\frac vL\right)
=\frac{\epsilon(\eta,v)}L n \partial_{v_2}f(v/L)+O\left(\left(\frac{k(\eta,v)}L\right)^2\right).
\end{eqnarray}
Altogether we get
\begin{eqnarray}
  \begin{aligned}
 \partial_t \mathbb E\left(I^{(L)}_f(\eta({L^2 t}))\right)&=  L^2 \mathbb E[\mathcal L^I I_f](\eta({L^2 t}))=\frac12 \frac1{|V_L|}\sum_{v\in V_L}\frac{\mathbb E\hat h_{\eta(L^2 t)}(v)}L \partial^2_{v_1}f(v/L)
\\&+ \frac1{|V_L|}\sum_{v\in V_L}\partial_{v_2}f(v/L)
\mathbb E\left[r(\eta(L^2t),v)\frac{k(\eta(L^2t),v)-1}2 \right]\\& +O\left(\frac1{L |V_L|}\sum_{v\in V_L}\mathbb E \,r(\eta(L^2t),v)k^2(\eta(L^2t),v)\right)  +O(1/L).
  \end{aligned}
\end{eqnarray}
The first sum in the r.h.s. tends to 
\begin{eqnarray}
  \frac12\int_V \hat \psi(v,t) \partial^2_{v_1}f(v) dv: 
\end{eqnarray}
just write the height $h_\eta(v)$ as the height at a point $v_0$ on
the boundary of $V_L$ plus the sum of interface gradients from $v_0$ to
$v$, then use the hypothesis of local equilibrium to replace the
average interface gradients in terms of derivatives of $\hat \psi$.

As for the second term, we claim that under the assumption of local
equilibrium we can replace
\begin{eqnarray}
  \label{eq:qqqq}
\frac1{|V_L|}\sum_{v\in V_L}\partial_{v_2}f(v/L)\mathbb E\left[r(\eta(L^2t),v)\frac{k(\eta(L^2t),v)-1}2 \right] 
\end{eqnarray}
with
\begin{eqnarray}
  \label{eq:repl2}
\int_V \partial_{v_2}f(v)\left(1+\frac12 \partial_{v_2}\hat\psi(v,t)\right)  
\pi_{\rho(\nabla\hat\psi(v,t))}\left[1_{X(0,0)}1_{n(b^{(0,0)})\ne0}\frac{|n(b^{(0,0)})|-1}2\right] dv.
\end{eqnarray}
Let us see why.  Let $v_L\in V_L$ be such that $v_L/L\to \bar v$ in
the interior of $V$ and let $V_{\epsilon,L}$ be a ball of radius
$\epsilon L$ around $v_L$. If $\Sigma_{\epsilon,L}$ is the portion of
the interface $\Sigma$ whose $\Pi_{110}$ projection is
$V_{\epsilon,L}$, let $U_{\epsilon,L}$ be its $\Pi_{111}$ projection.
Then we have  (with $\eta\equiv \eta(L^2 t)$)
\begin{multline}
  \frac1{|V_{\epsilon,L}|} \sum_{v\in V_{\epsilon,L}}
  r(\eta,v)\frac{k(\eta,v)-1}2\\=
  \frac{|U_{\epsilon,L}|}{|V_{\epsilon,L}|}\times
  \frac1{|U_{\epsilon,L}|}\sum_{u\in
    U_{\epsilon,L}}1_{X(u)}1_{n(b^{(u)})\ne
    (u_1+u_2)/2}\frac{|n(b^{(u)})-(u_1+u_2)/2|-1}2.
\label{eq:furbata}
\end{multline}
Indeed, the event $X(u)\cap \{n(b^{(u)})\ne (u_1+u_2)/2\}$ is the event that the edge $u,u+(1,1)$ is the common side of a lozenge of type $1$ and one of type $2$, see Definition \ref{def:Xu}.
In turn, this condition says that the point $u$ corresponds to one of the values of $v\in V_L$ with $r(\eta,v)\ne 0$, see Figure \ref{fig:ennesima}.

Since $|U_{\epsilon,L}|/|V_{\epsilon,L}|$ equals the ratio between the
total number of faces of $\Sigma_{\epsilon,L}$ and the number of its
non-horizontal faces, by \eqref{eq:15} and \eqref{eq:sbizzarro2} this ratio converges to
$[1+\frac12 \partial_{v_2}\hat\psi(\bar v,t)]$. As for the second term in
\eqref{eq:furbata}, by the assumption of local equilibrium and the
translation invariance of the Gibbs measures $\pi_\rho$, its average
converges to the average in \eqref{eq:repl2}.

Similarly to \eqref{eq:repl2}, one can replace 
\begin{eqnarray}
\label{eq:ggg}
  \frac1{L|V_L|}\sum_{v\in V_L}\mathbb E \,r(\eta(L^2t),v)k^2(\eta(L^2t),v)
\end{eqnarray}
with 
\begin{eqnarray}
\label{eq:gg}
 \frac1L \int_V\left(1+\frac12 \partial_{v_2}\hat\psi(v,t)\right)\pi_{\rho(\nabla\hat\psi(v,t))}\left[1_{X(0,0)}1_{n(b^{(0,0)})\ne0}\,n(b^{(0,0)})^2\right]dv=O(1/L)
\end{eqnarray}
(the random variable $n(b^{(0,0)})$ has moments of all orders under $\pi_\rho$).
Altogether, \eqref{eq:hypo2} follows.
\begin{Remark}
  \label{rem:tec}
  To justify that \eqref{eq:ggg} is $o(1)$ it would be clearly
  sufficient to prove that the second moment of $k(\eta,v)$ during the
  evolution remains $o(L)$. Under the equilibrium Gibbs
  measures $\pi_\rho$ the random variable $k(\eta,v)$ is known to have  exponential
  moments of all orders, so the situation looks promising.

  However, justifying the replacement of \eqref{eq:qqqq} with
  \eqref{eq:repl2} looks definitely harder: following the $H^{-1}$
  method of \cite{FuSpo}, it appears that one needs some form of uniform integrability 
  bound on $k(\eta,v)$, uniformly in $L$. In the forthcoming work \cite{cf:hydro}, in the case of periodic boundary conditions we manage to bypass  this difficulty.

\end{Remark}

\section{Some equilibrium computations}

\label{sec:K}

\subsection{Proof of Proposition  \ref{lemma:piX}}

\label{sec:piX}
We begin with proving Eq. \eqref{eq:piX}.
From the last remark in Definition \ref{def:Xu} we see that
the event $X(0,0)\cap \{n(b^{(0,0)})\ne0\}$ is equivalently
the event that the edge $e$ joining vertices $(0,0) $ and $(1,1)$ of $\mathcal T$ is the common side of two lozenges, one of type $1$ and one of type $2$.
From \eqref{eq:piro} we see that (see Figure \ref{fig:1} for notations) 
\begin{multline}
  \label{eq:pis1}
  \pi_\rho({X(0,0)}\cap \{n(b^{(0,0)})\ne0\})\\=2 k_1 k_2 [K^{-1}(w_1,b_1)K^{-1}(w_2,b_2)-K^{-1}(w_1,b_2)K^{-1}(w_2,b_1)],
\end{multline}
where the factor $2$ comes from the choice of whether the type-$1$ lozenge is to the right or to the left of the edge $e$.
From \eqref{eq:K1i},  we see that $k_1 K^{-1}(w_1,b_1)=\rho_1,k_2 K^{-1}(w_2,b_2)=\rho_2,k_3 K^{-1}(w_1,b_2)=(1-\rho_1-\rho_2)$. 
As a consequence, 
\begin{multline}
  \label{eq:pis2}
    \pi_\rho(X(0,0)\cap\{ n(b^{(u)})\ne0\})=2(\rho_1 \rho_2 + (1-\rho_1-\rho_2)\times \Xi),\\
 \Xi=- \frac{k_1k_2}{k_3}\frac1{(2\pi i)^2}
 \oint dz \oint dw \frac1{k_3+k_1 z+k_2 w}.
\end{multline}
Thanks to \cite[Theorem 2.7]{FerrariSunil}\footnote{In the notations of \cite{FerrariSunil}, $k_1=b,k_2=c$ and $k_1=a$}, we know that $\Xi$ equals exactly $V(\rho)$.

Next we prove \eqref{eq:orr}.
For lightness of notation, we write $n, n^+,n^-,X,|I|$ instead of
$n(b^{(0,0)}),n^+(b^{(0,0)}),n^-(b^{(0,0)}), X(0,0), |I(b^{(0,0)})|$.
\begin{Remark}
\label{rem:uniformi}
Since $\pi_\rho$ is translation-invariant and is locally the uniform measure on lozenge tilings, we
have:
\begin{enumerate}
\item [(U1)] conditionally on the event 
$X$ and on the values of $n^+,n^-$, the random variable $n$ is uniformly distributed in $\{n^-,\dots, n^+\}$.
\item [(U2)] conditionally on the event 
$X$ and on $n^+-n^-$, the random variable $n^+$ is uniformly distributed on $\{0,\dots, n^+-n^-\}$. 
\end{enumerate}
  
\end{Remark}

We want to compute
\begin{eqnarray}
\label{eq:K'}
  \pi_\rho((|n|-1)1_{X}1_{n\ne0})= \pi_\rho(|n|1_X 1_{n\ne0})-\pi_\rho(1_{X}1_{n\ne0}).
\end{eqnarray}
We have 
\begin{eqnarray} 
\label{eq:1/3}
\pi_\rho(|n|1_X 1_{n\ne0})=\pi_\rho(|n| 1_{X})=\frac13\pi_\rho\left(\frac{|I|^2-1}{|I|}1_{X}\right).
\end{eqnarray}
Indeed, from (U1)-(U2) we see that, conditionally on
$|I|:=n^+-n^-+1$, the two random variables $n^+$ and $n^+-n$ are 
independent and uniformly distributed 
in
$\{0,\dots,|I|-1\}$. We have already mentioned in Section \ref{sec:GK} that  if $u_1,u_2$ are independent and uniform in $\{0,\dots,|I|-1\}$ then
\[
E(|u_1-u_2|)=\frac13 \frac{|I|^2-1}{|I|}.
\]
Using (U1), the last expression in \eqref{eq:1/3} is also
\begin{eqnarray}
\frac13\pi_\rho((|I|+1)1_{X}1_{n\ne0})  
=2V(\rho)
\end{eqnarray}
where in the last step we used Theorem \ref{prop:K}.

On the other hand, $
  \pi_\rho(1_{X}1_{n\ne0})$ was computed in  \eqref{eq:piX}. Altogether \eqref{eq:K'} equals
\begin{eqnarray}
  2[-\rho_1 \rho_2+(\rho_1+\rho_2)V(\rho)]
\end{eqnarray}
and \eqref{eq:orr} is proven.

\subsection{Proof of Theorem  \ref{prop:K}}
Again, we write $n, n^+,n^-,X$ instead of
$n(b^{(0,0)}),n^+(b^{(0,0)})$, $n^-(b^{(0,0)}), X(0,0)$.
Then, Eq. \eqref{eq:1} together with definition \eqref{eq:CF} gives 
  \begin{eqnarray}
    \label{eq:FS}
 \Delta:=   \pi_\rho\left[(n^++1)1_{X}1_{n<0}\right]=\frac1\pi \frac{\sin(\pi \rho_1)\sin(\pi \rho_2)}{\sin(\pi(1-\rho_1-\rho_2))}
  \end{eqnarray}
and we need to prove that this equals
\begin{eqnarray}
\label{58}
  \frac1{6}\pi_\rho [(n^+-n^-+2)\,1_X 1_{n\ne 0}]
\end{eqnarray}
which is twice the r.h.s. of \eqref{eq:mu}.

The Gibbs measures $\pi_\rho$ are invariant by reflection through any
vertex of $\mathcal T$ (see Remark 3.3 in \cite{Toninelli2+1}).  As a
consequence of this symmetry, conditionally on $X$, the random
variables $n^+\ge0$ and $-n^-\ge0$ have the same law and
\begin{eqnarray}
\Delta=  \pi_\rho\left[(
1-n^-)1_{X}1_{n>0}
\right].
\end{eqnarray}
We deduce from property (U1) above that 
\begin{eqnarray}
  2\Delta=\pi_\rho\left[(1-n^-)\frac{n^+}{n^+-n^-+1}+
(1+n^+)\frac{|n^-|}{n^+-n^-+1};X
\right].
\end{eqnarray}
In fact, by first conditioning on the event $X$ and then on  the value of $n^+,n^-$ and using that $n$ is conditionally uniform on $\{n^-,\dots,n^+\}$ one has (with $\pi_\rho^X$ the measure $\pi_\rho$ conditioned to the event $X$)
\begin{multline}
  \pi_\rho[(n^++1)1_X1_{n<0}]=\pi^X_\rho[(n^++1)1_{n<0}]\pi_\rho(X)\\=\pi^X_\rho\left[(n^++1)
\pi^X_\rho(1_{n<0}|n^+,n^-)
\right]\pi_\rho(X)=\pi^X_\rho\left[(n^++1)
\frac{|n^-|}{n^+-n^-+1}
\right]\pi_\rho(X)\\=\pi_\rho\left[
(n^++1)
\frac{|n^-|}{n^+-n^-+1}
\right],
\end{multline}
and similarly for $ \pi_\rho[(1-n^-)1_X1_{n>0}]$.
Using property (U2) (together with $E U=k/2, E U^2=(1/6)k(1+2k)$ for
a uniform random variable $U$ on $\{0,\dots,k\}$) we have
\begin{eqnarray}
  \pi_\rho[n^+ n^-|X,n^+-n^-]=-\frac16(n^+-n^-)^2+\frac16(n^+-n^-).
\end{eqnarray}
As a consequence, 
\begin{gather}
  \label{eq:deltat}
  \Delta= \frac1{6}\pi_\rho\left[\frac{(n^+-n^-)^2+2(n^+-n^-)}{n^+-n^-+1}
;X
\right]= \frac1{6} \pi_\rho\left[(n^+-n^-+2)\frac{n^+-n^-}{n^+-n^-+1}; X\right].
\end{gather}
On the other hand, from (U2),   this equals \eqref{58}.

\section*{Acknowledgments}
We would like to thank the anonymous referee for a careful reading.
We are grateful to  Pietro Caputo and Fabio Martinelli for
innumerable enlightening discussions.   F.  T.  was  partially
funded by Marie Curie IEF Action DMCP ``Dimers, Markov chains and Critical Phenomena'', grant
agreement n.  621894 and by the CNRS PICS grant “Interfaces al\'eatoires discr\`etes et dynamiques
de Glauber”. B. Laslier was supported by the Engineering and Physical Sciences
Research Council under grant EP/103372X/1.

\end{document}